\documentclass[12pt, a4paper]{article}
\title{Some applications of the multiplicative tensor product of matrix factorizations}
\author{Yves Baudelaire Fomatati\\
\small Department of Mathematics and Statistics, University of Ottawa,\\ \small Ottawa, Ontario, Canada K1N 6N5.\\ \small yfomatat@uottawa.ca.}

\date{}

\usepackage[square,comma,numbers,sort&compress]{natbib}
\usepackage[centertags]{amsmath}
\usepackage{latexsym}
\usepackage{leqno}
\usepackage{amsfonts}
\usepackage{amsmath,amssymb}
\usepackage{enumerate}

\usepackage{amsthm}
\usepackage[utf8,latin9]{inputenc}
\usepackage{xcolor}

\theoremstyle{plain}
\newtheorem{remark}{Remark}[section]
\theoremstyle{plain}

\theoremstyle{plain}
\newtheorem{proposition}{Proposition}[section]
\theoremstyle{plain}
\newtheorem{theorem}{Theorem}[section]
\theoremstyle{plain}
\newtheorem{definition}{Definition}[section]
\theoremstyle{plain}

\theoremstyle{plain}

\theoremstyle{plain}
\newtheorem{notation}{Notations}[section]
\newtheorem{example}{Example}[section]

\theoremstyle{plain}

\usepackage[all, 2cell]{xy}
\usepackage{txfonts}

\begin{document}
\maketitle
\begin{quote}
  \textbf{Abstract}
\end{quote}
The notion of semi-unital semi-monoidal category was defined a couple of years ago using the so called "Takahashi tensor product" and so far, the only example of it in the literature is complex. In this paper, we use the recently defined "multiplicative tensor product of matrix factorizations" to give a simple example of this a notion. In fact,
if $MF(1)$ denotes the category of matrix factorizations of the constant power series $1$, we define the concept of one-step connected category and prove that there is a one-step connected subcategory of $(MF(1),\widetilde{\otimes})$ which is semi-unital semi-monoidal.
We also define the concept of right pseudo-monoidal category which generalizes the notion of monoidal category and we prove that $(MF(1),\widetilde{\otimes})$ is an example of this concept. 

\section{Introduction}

Eisenbud introduced the concept of matrix factorization (cf. \cite{eisenbud1980homological}) in 1980. His research results show how to use matrices to factorize all polynomials. For example the irreducible polynomial $f(x)=x^{2}+y^{2}$ over $\mathbb{R}[x,y]$ can be factored as follows: $$\begin{bmatrix}
    x  &  -y      \\
    y  &  x
\end{bmatrix}
\begin{bmatrix}
    x  &  y      \\
    -y  & x
\end{bmatrix}
= (x^{2} + y^{2})\begin{bmatrix}
    1  &  0      \\
    0  &  1
\end{bmatrix}
=fI_{2} $$
Thus, we say that
$
(\begin{bmatrix}
    x  &  -y      \\
    y  &  x
\end{bmatrix},
\begin{bmatrix}
    x  &  y      \\
    -y  & x
\end{bmatrix})
$
 is a $2\times 2$ matrix factorization of $f$.\\
 In a sense, this notion of factorizing polynomials using matrices can be seen as a generalization of the classical notion of polynomial factorization where a polynomial $p(x)=q(x)r(x)$ can be considered as the product of two $1\times 1$ matrices.\\
Matrix factorizations play an important role in many areas of pure mathematics
and physics. The notion of matrix factorizations is one of the key tools used in representation theory of hypersurface rings. It is a classical tool in the study of hypersurface singularity algebras (cf. \cite{eisenbud1980homological}).
 One of the discoveries of Eisenbud is that matrix factorizations of $f\in K[[x]]$ are closely related to the homological properties of modules over quotient rings $K[[x]]/(f)$.
\\ Let $K$ be a field and $K[[x]]$ be the formal power series ring in the variables $x=x_{1},\cdots,x_{r}$ and $K[[y]]$ be the formal power series ring in the variables $y=y_{1},\cdots,y_{s}$. Let $f\in K[[x]]$ and $g\in K[[y]]$ be nonzero noninvertible\footnote{Yoshino \cite{yoshino1998tensor} requires an element $f\in K[[x]]$ to be nonzero noninvertible because if $f=0$ then $K[[x]]/(f)=K[[x]]$ and if $f$ is a unit, then $K[[x]]/(f)=$K[[x]]/K[[x]]$=\{1\}$. But in this paper we will not bother about such restrictions because we will not deal with the homological methods used in \cite{yoshino1998tensor}.} elements.
In 1998, Yoshino constructed a tensor product denoted $\widehat{\otimes}$ which is such that if $X$ is a matrix factorization of $f\in K[[x]]$ and $Y$ is a matrix factorization of $g\in K[[y]]$, then $X\widehat{\otimes}Y$ is a matrix factorization of $f+g\in K[[x,y]]$. In 2002 and 2003,  Kapustin and Li in their papers \cite{kapustin2003topological} and \cite{kapustin2004d}, used matrix factorizations in string theory to study boundary conditions for strings in Landau-Ginzburg models. In 2012, Carqueville and Murfet in their paper \cite{carqueville2016adjunctions}, briefly presented the construction of the bicategory $\mathcal{LG}_{K}$ of Landau-Ginzburg models
whose $1$-cells are matrix factorizations. In 2013, the geometry of the category of matrix factorizations was studied in Yu's Ph.D. dissertation \cite{yu2013geometric}. In 2014, Camacho \cite{camacho2015matrix} in her PhD dissertation recalled the notion of graded matrix factorizations with special emphasis on $\mathbb{C}-$graded matrix factorizations.
\\
In 2016, Crisler and Diveris \cite{crisler2016matrix} examined matrix factorizations of polynomials in the ring $\mathbb{R}[x_{1},\cdots,x_{n}]$, using only techniques from elementary linear algebra. They focused mostly on factorizations of sums of squares of polynomials. They improved the standard method for factoring polynomials for this class of polynomials.
More recently in 2019, the author in his Ph.D. dissertation \cite{fomatati2019multiplicative} defined the multiplicative tensor product of matrix factorizations and found a variant of this product \cite{fomatati2021tensor} in 2020. These were then used to further improve the standard method for factoring a large class of polynomials. In \cite{fomatati2021necessary}, properties of matrix factorizations are used to find a necessary condition to obtain a Morita Context in the bicategory of Landau-Ginzburg models.
\\
 In this paper, we use the recently defined multiplicative tensor product of matrix factorizations \cite{fomatati2019multiplicative} to give a simple example of a semi-unital semi-monoidal category. This is a notion that was defined recently in \cite{abuhlail2013semiunital} using the so called 'Takahashi tensor product' and it required a complex set-up. We will construct an easy-to-understand example with a relatively small amount of set-up. Moreover, we will also use the multiplicative tensor product of matrix factorization to define and give an example of the concept of right pseudo-monoidal category which generalizes the notion of monoidal category. 
 \\

 \textbf{Significance of the notion of semi-unit:} \\
We first recall some definitions.\\
A \textit{semi-ring} is roughly speaking, a ring not necessarily with subtraction. The first natural example of a semi-ring is the set $\mathbb{N}$ of non-negative integers.\\
A \textit{semi-module} is roughly speaking a module not necessarily with subtraction. The category of Abelian groups is nothing but the category of modules over $Z$; similarly, the category of commutative monoids is nothing but the category of semi-modules over $\mathbb{N}$\\
Semi-rings were studied by many algebraists beginning with Dedekind. They have significant applications in several areas, for instance  Automata  Theory and  Optimization  Theory (see \cite{golan1999semirings} for applications).\\
The theory of semi-modules over semi-rings was developed by many authors including Takahashi \cite{takahashi1982bordism}.
 In 2008, Jawad used the so called Takahashi's tensor-like product $\boxtimes$ of semi-modules over an associative semi-ring $A$ \cite{takahashi1982bordism}, to introduce notions of semi-unital semi-rings and semi-counital semi-corings (cf. \cite{abuhlail2013semiunital}). However, these could not be realized as monoids (comonoids) in the category $_{A}S_{A}$ of $(A,A)$-bisemi-modules. This is mainly due to the fact that the category $(_{A}S_{A}, \boxtimes, A)$ is not monoidal in general. Motivated by the desire to fix this problem, Jawad \cite{abuhlail2013semiunital} introduced and investigated a notion of semi-unital semi-monoidal categories with prototype $(_{A}S_{A}, \boxtimes, A)$ and investigated semi-monoids (semi-comonoids) in such categories as well as their categories of semi-modules (semi-comodules). He realized that although the base semi-algebra $A$ is not a unit in $_{A}S_{A}$, $A$ nevertheless has properties of what he called a \textit{semi-unit}. This motivated the introduction of a more generalized notion of monads (comonads) in arbitrary categories (for more on this, see \cite{abuhlail2013semiunital}).
\begin{example}
An example of semi-unital semi-monoidal category as given by Jawad in \cite{abuhlail2013semiunital} is the category of bisemi-modules over a semi-algebra $A$ with the Takahashi tensor product $(_{A}S_{A}, \boxtimes, A)$. That is the only example we found in the literature. unfortunately, it requires a great amount of set-up and so we refer the reader to theorem 5.11 of \cite{abuhlail2013semiunital}. As earlier mentioned, a (less involved) example of semi-unital semi-monoidal category will be given in this paper (cf. theorem \ref{Thm: T is a semi-unital semi-monoidal catg}) using the recently defined multipicative tensor product of matrix factorizations \cite{fomatati2019multiplicative}.
\end{example}

In the next section, definitions of special classes of categories are recalled. In section 3, the notion of tensor products of matrix factorizations in also recalled. A comparison of the tensor product of matrix factorizations and its multiplicative counterpart is presented in section 4. The category of matrix factorizations of the constant power series $1$ is studied under this section. Moreover, a simple example of a semi-unital semi-monoidal category using the multiplicative tensor product is presented. We wrap up this section with the definition of the notion of right pseudo-monoidal category. 

\section{Special classes of categories}
Here, we recall the definitions of some special types of categories.
\begin{definition}\cite{mac2013categories} \label{defn monoidal category}
A \textbf{monoidal category} $\mathcal{C}=<\mathcal{C},\square, e, \alpha, \lambda, \rho>$ is a category $\mathcal{C}$, a bifunctor $\square : \mathcal{C} \times \mathcal{C}\rightarrow \mathcal{C}$, an object $e\,\in \,\mathcal{C}$, and
three natural isomorphisms $\alpha,\,\lambda,\,and\,\rho$;\\
such that:\\
$\bullet$  $\alpha=\alpha_{a,b,c}: a\square (b\square c) \cong (a\square b)\square c$ is natural for all $a,b,c\,\in \,\mathcal{C}$ and the pentagonal diagram

$\xymatrix @ R=0.4in @ C=1.0in{
a\square (b\square (c\square d)) \ar[rr]^{1_{a}\square \alpha}\ar[d]_\alpha &&
a\square ((b\square c)\square d) \ar[d]^\alpha &\\
(a\square b)\square (c\square d) \ar[rd]_\alpha &&
(a\square(b\square c))\square d \ar[ld]^{\alpha\square 1_{d}} &\\
&((a\square b)\square c)\square d} $

commutes for all $a,b,c,d\,\in \,\mathcal{C}$.\\
($\alpha$ is also called associator (p.11 \cite{carqueville2016adjunctions})).\\
$\bullet$ $\lambda$ and $\rho$ are natural. On p.10 of \cite{carqueville2016adjunctions}, $\lambda$ and $\rho$ are also called left and right unit actions (or unitors).
\\
$\lambda_{a}: e\square a \cong a$, $\rho_{a}: a\square e \cong a$,
for all objects $a\,\in \mathcal{C}$. \\Moreover, the triangular diagram
\[
\xymatrix{
a\square (e\square b)\ar[rr]^-\alpha \ar[rd]_{1_{a}\square \lambda_{b}} &&
(a\square e)\square b \ar[ld]^{\rho_{a}\square 1_{b}}\\
&a\square b}
\]
commutes for all $a,b\,\in \,\mathcal{C}$ and $\lambda_{e}=\rho_{e}:e\square e \rightarrow e$.
\end{definition}
\begin{definition}\cite{mac2013categories} \label{defn symm monoidal category}
A \textbf{symmetric monoidal category} is a monoidal category together with a symmetry. A symmetry $\gamma$
for a monoidal category $\mathcal{C}=<\mathcal{C},\square, e, \alpha, \lambda, \rho>$ is a natural isomorphism $\gamma = \gamma_{ab}: a\square b \rightarrow b\square a$ such that the following three diagrams commute for all $a,b,c \,\in \,\mathcal{C}$:
\[
\xymatrix{
a\square b\ar[rr]^-\gamma \ar[rd]_{1} &&
b\square a \ar[ld]^{\gamma_{ba}}\\
&a\square b}
\]
\[
\xymatrix{
e\square b\ar[rr]^-\gamma \ar[rd]_{\lambda_{b}} &&
b\square e \ar[ld]^{ \rho_{b}}\\
& b}
\]

$\xymatrix @ R=0.4in @ C=1.0in{
a\square(b\square c) \ar[rr]^{1_{a}\square \gamma}\ar[d]_\alpha &&
a\square (c\square b) \ar[d]^\alpha &\\
(a\square b)\square c \ar[d]_\gamma &&
(a\square c)\square b \ar[d]^{\gamma\square 1_{b}} &\\
c\square (a\square b) \ar[rr]^{\alpha}&&(c\square a)\square b}$
\end{definition}
An example of a symmetric monoidal category is the category of vector spaces over some fixed field $K$, using the ordinary tensor product of vector spaces.

\begin{definition} \label{defn semi-monoidal category} \cite{kock2008elementary}
  A \textbf{semi-monoidal category} $\mathcal{C}=<\mathcal{C},\square, \alpha>$ is a category $\mathcal{C}$, a bifunctor $\square : \mathcal{C} \times \mathcal{C}\rightarrow \mathcal{C}$ and
a natural transformation $\alpha$, satisfying the following condition:\\
$\bullet$ $\alpha$ is a natural isomorphism with components $\alpha_{a,b,c}:(a\square b)\square c \rightarrow a\square (b\square c)$ such that the following pentagonal diagram

$\xymatrix @ R=0.4in @ C=0.80in{
a\square (b\square (c\square d)) \ar[rr]^{1_{a}\square \alpha}\ar[d]_\alpha &&
a\square ((b\square c)\square d) \ar[d]^\alpha &\\
(a\square b)\square (c\square d) \ar[rd]_\alpha &&
(a\square(b\square c))\square d \ar[ld]^{\alpha \square 1_{d}} &\\
&((a\square b)\square c)\square d} $

commutes for all $a,b,c,d\,\in \,\mathcal{C}$.
($\alpha$ is also called associator (p.11 \cite{carqueville2016adjunctions})).
\end{definition}

 \begin{definition} (Defn 3.1 \cite{selinger2008idempotents})
 Given a category $\mathcal{C}$ and an object $A$ of $\mathcal{C}$, an \textbf{idempotent} of $\mathcal{C}$ is an endomorphism
 $e:A\rightarrow A$ with $e \circ e = e$.
An idempotent $e:A\rightarrow A$ is said to \textbf{split} if there is an object $B$ and morphisms $r: A \rightarrow B$, $s: B \rightarrow A$ such that $e = s\circ r$ and $id_{B} = r\circ s$.
\end{definition}

\begin{remark} (cf. Remark 3.4 \cite{selinger2008idempotents})
The splitting of an idempotent is a special case of a categorical limit and colimit.  More precisely, if $e:A\rightarrow A$ is an idempotent, then $r: A \rightarrow B$, $s: B \rightarrow A$ is a splitting of $e$ if and only if $r$ is a colimit and $s$ is a limit of the diagram $e:A\rightarrow A$.
\end{remark}

\begin{definition} \cite{borceux1986cauchy}
 An ordinary category is \textbf{idempotent complete}, aka \textbf{Karoubi complete} or \textbf{Cauchy complete}, if every idempotent splits.
 \end{definition}

The phrase "ordinary category" as opposed to the phrase "Higher category" is used here. Higher category theory is the generalization of category theory to a context where there are not only morphisms between objects, but generally $n$-morphisms between $(n-1)$-morphisms, for all $n \in \mathbb{N}$.


\begin{definition} \label{defn semi-unital semimonoidal catg} \cite{abuhlail2013semiunital}
Let $(\mathcal{C},\Box)$ be a semi-monoidal category with natural isomorphism $\alpha_{X,Y,Z}: (X\Box Y)\Box Z\rightarrow X\Box (Y\Box Z)$ for all $X,Y,Z \in \mathcal{C}$. Let $\mathbb{I}$ stand for the identity endofunctor on any given category. We say that $\mathbf{I}\in \mathcal{C}$ is a \textbf{semi-unit} if the following conditions hold:
\begin{enumerate}
  \item There is a natural transformation $\omega: \mathbb{I} \rightarrow (\mathbf{I}\Box -)$;
  \item There exists an isomorphism of functors $\mathbf{I}\Box - \cong -\Box \mathbf{I}$, i.e., there is a natural isomorphism $l_{X}: \mathbf{I}\Box X \cong X\Box \mathbf{I}$ in $\mathcal{C}$ with inverse $q_{X}$, for each object $X$ of $\mathcal{C}$, such that $l_{\mathbf{I}}=q_{\mathbf{I}}$ and the following diagrams are commutative for all $X,Y \in \mathcal{C}$:

 $$\xymatrix@ R=0.6in @ C=.75in{(\mathbf{I}\Box a)\Box b\ar[r]^{\alpha_{\mathbf{I},a,b}} \ar[d]_{l_{a}\Box b} &
\mathbf{I}\Box (a\Box b) \ar[r]^{l_{a\Box b}} & (a\Box b)\Box \mathbf{I}\ar[d]^{\alpha_{a,b,\mathbf{I}}\,\,\,\,\,\,\,\,\,\,\,\,\,\,\,\,(1)}\\
(a\Box \mathbf{I})\Box b \ar[r]^{\alpha_{a,\mathbf{I},b}} & a\Box (\mathbf{I}\Box b)\ar[r]^{a \Box l_{b}} & a\Box (b\Box \mathbf{I})}$$

\[
\xymatrix{
a\Box b\ar[rr]^{\omega_{a}\Box b} \ar[rd]_{\omega_{a\Box b}} &&
(\mathbf{I}\Box a) \Box b \ar[ld]^{\cong \,\,\,\,\,\,\,\,\,\,\,\,\,\,\,\,(2)}\\
&\mathbf{I}\Box (a \Box b)}
\]
\[
\xymatrix{
a\Box b\ar[rr]^{a \Box\omega_{b}} \ar[rd]_{\omega_{a\Box b}} &&
a\Box (\mathbf{I} \Box b) \ar[ld]^{\cong \,\,\,\,\,\,\,\,\,\,\,\,\,\,\,\,(3)}\\
&\mathbf{I}\Box (a \Box b)}
\]
\end{enumerate}
A \textbf{semi-unital semi-monoidal category} is a semi-monoidal category with a semi-unit.
\end{definition}
\begin{remark}\cite{abuhlail2013semiunital}
  If $X\cong \mathbf{I}\Box X (\cong X\Box \mathbf{I})$, then we say that $X$ is \textbf{firm} and set $\lambda_{X}:= \omega_{X}^{-1}: \mathbf{I}\Box X \rightarrow X$. \\
  If $\mathbf{I}$ is firm (also called a \textbf{pseudo-idempotent}) and $\omega_{\mathbf{I}}^{-1}\Box \mathbf{I}=\mathbf{I}\Box \omega_{\mathbf{I}}^{-1}$, then $\mathbf{I}$ is an \textbf{idempotent}.\\
  A semi-monoidal category becomes a monoidal category if it has a unit, i.e., $\mathbf{I}$ is such that $\lambda_{\mathbf{I}}=\kappa_{\mathbf{I}}$, $\lambda_{X\Box Y}=\lambda_{X}\Box Y$
  and $\kappa_{X\Box Y}=X\Box \kappa_{Y}$ for all $X,Y\in \mathcal{C}$, where $\kappa_{X}=\omega_{X}^{-1}\circ q_{X}: X\Box \mathbf{I} \cong \mathbf{I}\Box X \rightarrow X$.
\end{remark}
\begin{remark} (cf. remark 3.3 of \cite{abuhlail2013semiunital}, \cite{kock2008elementary}) \label{remark on units}
  Kock in \cite{kock2008elementary} called an object $I$ a \textbf{Saavedra unit} or \textbf{reduced unit} just in case it is pseudo-idempotent and cancellable in the sense that the endofunctors $I\Box -$ and $- \Box I$  are full and faithful (equivalently, $I$ is idempotent and the endofunctors $I\Box -$ and $- \Box I$ are equivalences of categories). Kock also showed that $I$ is a unit if and only if it is a Saavedra unit. \\
  Since every unit is a semi-unit, the notion of semi-unital semi-monoidal categories generalizes the classical notion of monoidal categories.
\end{remark}

There is also a notion of skew-monoidal category (\cite{szlachanyi2012skew}) defined as follows:
\begin{definition} (cf. \cite{szlachanyi2012skew}) \label{defn: right monoidal catg}
A \textbf{right-monoidal category} $(\mathcal{M},\ast,e,\alpha,\gamma,\rho)$ consists of a category $\mathcal{M}$, a functor $(-)\ast(-): \mathcal{M}\times \mathcal{M} \rightarrow \mathcal{M}$, an object $e$ of $\mathcal{M}$ and natural transformations\\
$\alpha_{L,M,N}:L\ast (M\ast N) \rightarrow (L\ast M)\ast N$ \\
$\gamma_{M}: M \rightarrow e\ast M$ \\
$\rho_{M}: M\ast e \rightarrow M$

subject to the following axioms: For all objects $K,L,M,N$
\begin{enumerate}
  \item $(\alpha_{K,L,M}\ast N)\circ \alpha_{K,L\ast M,N}\circ (K\ast \alpha_{L,M,N})= \alpha_{K\ast L,M,N}\circ \alpha_{K,L,M\ast N}$
  \item $\alpha_{e,M,N}\circ \gamma_{M\ast N}= \gamma_{M}\ast N$
  \item $\rho_{M\ast N}\circ \alpha_{M,N,e}= M\ast \rho_{N}$
  \item $(\rho_{M}\ast N)\circ\alpha_{M,e,N}\circ (M\ast \gamma_{N})=id_{M\ast N}$
  \item $\rho_{e}\circ \gamma_{e}= id_{e}$.
\end{enumerate}
\end{definition}
\begin{remark} (cf. \cite{szlachanyi2012skew})
  If $\mathcal{M}$ is replaced with $\mathcal{M}^{op}$, we obtain what is called a \textbf{left-monoidal category}.\\
If $\alpha$, $\gamma$ and $\rho$ are isomorphisms, we recover the notion of monoidal category. So a right-monoidal category is a generalization of a monoidal category.
\end{remark}

\section{Tensor products of matrix factorizations}
In this section, we recall the definitions of the Yoshino tensor product of matrix factorizations denoted $\widehat{\otimes}$. Next, we recall the definition of \textit{multiplicative tensor product} of matrix factorizations denoted $\widetilde{\otimes}$. \\
Under this section, unless otherwise stated, $R=K[[x]]$ and $S=K[[y]]$ where $x=x_{1},...,x_{r}$ and $y=y_{1},...,y_{s}$.

\subsection{Yoshino's tensor product of matrix factorization} \label{subsec: Yoshino tens prod and variant}
Recall the following:
\begin{definition} \cite{yoshino1998tensor}  \label{defn matrix facto of polyn}   \\
An $n\times n$ \textbf{matrix factorization} of a power series $f\in \;R$ is a pair of $n$ $\times$ $n$ matrices $(\phi,\psi)$ such that
$\phi\psi=\psi\phi=fI_{n}$, where $I_{n}$ is the $n \times n$ identity matrix and the coefficients of $\phi$ and of $\psi$ are taken from $R$.
\end{definition}
Also recall ($\S 1$ of \cite{yoshino1998tensor}) the definition of
the category of matrix factorizations of a power series $f\in R=K[[x]]:=K[[x_{1},\cdots,x_{n}]]$ denoted by $MF(R,f)$ or $MF_{R}(f)$, (or even $MF(f)$ when there is no risk of confusion):\\
$\bullet$ The objects are the matrix factorizations of $f$.\\
$\bullet$ Given two matrix factorizations of $f$; $(\phi_{1},\psi_{1})$ and $(\phi_{2},\psi_{2})$ respectively of sizes $n_{1}$ and $n_{2}$, a morphism from $(\phi_{1},\psi_{1})$ to $(\phi_{2},\psi_{2})$ is a pair of matrices $(\alpha,\beta)$ each of size $n_{2}\times n_{1}$ which makes the following diagram commute \cite{yoshino1998tensor}:
$$\xymatrix@ R=0.6in @ C=.75in{K[[x]]^{n_{1}} \ar[r]^{\psi_{1}} \ar[d]_{\alpha} &
K[[x]]^{n_{1}} \ar[d]^{\beta} \ar[r]^{\phi_{1}} & K[[x]]^{n_{1}}\ar[d]^{\alpha}\\
K[[x]]^{n_{2}} \ar[r]^{\psi_{2}} & K[[x]]^{n_{2}}\ar[r]^{\phi_{2}} & K[[x]]^{n_{2}}}$$
That is,
$$\begin{cases}
 \alpha\phi_{1}=\phi_{2}\beta  \\
 \psi_{2}\alpha= \beta\psi_{1}
\end{cases}$$
More details on this category are found in chapter 2 of \cite{fomatati2019multiplicative}.
\\
\begin{definition} \cite{yoshino1998tensor} \label{defn Yoshino tensor prodt}
Let $X=(\phi,\psi)$ be an $n\times n$ matrix factorization of $f\in R$  and $X'=(\phi',\psi')$ an $m\times m$ matrix factorization of $g\in S$. These matrices can be considered as matrices over $L=K[[x,y]]$ and the \textbf{tensor product} $X\widehat{\otimes} X'$ is given by\\
(\(
\begin{bmatrix}
    \phi\otimes 1_{m}  &  1_{n}\otimes \phi'      \\
   -1_{n}\otimes \psi'  &  \psi\otimes 1_{m}
\end{bmatrix}
,
\begin{bmatrix}
    \psi\otimes 1_{m}  &  -1_{n}\otimes \phi'      \\
    1_{n}\otimes \psi'  &  \phi\otimes 1_{m}
\end{bmatrix}
\))\\
where each component is an endomorphism on $L^{n}\otimes L^{m}$.
\end{definition}
It is easy to verify that $X\widehat{\otimes} X'$ is in fact an object of $MF_{L}(f+g)$ of size $2nm$.

\begin{remark}
When $n=1$, we get a $1$ $\times$ $1$ matrix factorization of $f$, i.e., $f=[g][h]$ which is simply a factorization of $f$ in the classical sense. But in case $f$ is not reducible, this is not interesting, that's why we will mostly consider $n > 1$.
\end{remark}
Variants of Yoshino's tensor product of matrix factorizations were constructed in \cite{fomatati2021tensor}. \\

\subsection{Multiplicative tensor product of matrix factorization}
In this subsection, we recall the definition of the multiplicative tensor product of matrix factorizations.
\\
First it is well known that if $A$ (resp. $B$) is an $m\times n$ (resp. $p\times q$) matrix, then their direct sum $A\oplus B = \begin{bmatrix}
    A     &       0      \\
    0     &       B
\end{bmatrix}$, where the $0$ in the first line is a $p\times q$ matrix and the one in the second line is an $m\times n$ matrix. \\
Finally, recall that if $A$ (resp. $B$) is an $m\times n$ (resp. $p\times q$) matrix, then their tensor product $A\otimes B$ is the matrix obtained by replacing each entry $a_{ij}$ of $A$ with the matrix $a_{ij}B$. So, $A\otimes B$ is a $mp\times nq$ matrix.

\begin{definition} \label{defn of the multiplicative tensor product} \cite{fomatati2019multiplicative}
Let $X=(\phi,\psi)$ be a matrix factorization of $f\in K[[x]]$ of size $n$ and let $X'=(\phi',\psi')$ be a matrix factorization of $g\in K[[y]]$ of size $m$. Thus, $\phi,\psi,\phi' \,and \,\psi'$ can be considered as matrices over $L=K[[x,y]]$ and the \textbf{multiplicative tensor product} $X\widetilde{\otimes} X'$ is given by \\\\
\[((\phi\otimes\phi')\oplus (\phi\otimes\phi'), (\psi\otimes\psi')\oplus (\psi\otimes\psi'))=(
\begin{bmatrix}
    \phi\otimes\phi'  &          0      \\
    0                  &\phi\otimes\phi'
\end{bmatrix},
\begin{bmatrix}
    \psi\otimes\psi'  &    0      \\
    0                 &  \psi\otimes\psi'
\end{bmatrix}
)\]
\\\\
where each component is an endomorphism on $L^{n}\otimes_{L} L^{m}$.
\end{definition}
\begin{remark}
One of the reasons for the "doubling" in the definition of the multiplicative tensor product of matrix factorizations is found in the proof of theorem \ref{Thm: T is a semi-unital semi-monoidal catg}. Observe that, in this proof, had we defined $X\widetilde{\otimes} X'$ as \[((\phi\otimes\phi'), (\psi\otimes\psi')),\] we would have had only one object in the category $\mathcal{T}$. Consequently, we would not have been able to construct another example of semi-unital semi-monoidal category.

\end{remark}

\begin{definition} \label{defn zeta tensor tilda X'} \cite{fomatati2019multiplicative}
For a morphism $\zeta=(\alpha, \beta): X_{1}=(\phi_{1},\psi_{1}) \rightarrow X_{2}=(\phi_{2},\psi_{2})$ in $MF(K[[x]],f)$ and for any $m\times m$ matrix factorization $X'=(\phi', \psi')$ in $MF(K[[y]],g)$, we define $\zeta \widetilde{\otimes} X'$
by
$$(
\begin{bmatrix}
  \alpha \otimes 1_{m}  &          0      \\
    0               &  \alpha \otimes 1_{m}
\end{bmatrix},
\begin{bmatrix}
    \beta \otimes 1_{m}   &    0      \\
    0                 &   \beta \otimes 1_{m}
\end{bmatrix}
)$$
\end{definition}

\begin{definition}\label{defn zeta is a funtor (right)} \cite{fomatati2019multiplicative}
For a morphism $\zeta'=(\alpha', \beta'): X_{1}'=(\phi_{1}',\psi_{1}') \rightarrow X_{2}'=(\phi_{2}',\psi_{2}')$ in $MF(K[[y]],g)$ and for any $n\times n$ matrix factorization $X=(\phi, \psi)$ in $MF(K[[x]],f)$, we define $X \widetilde{\otimes} \zeta'$
by
$$(
\begin{bmatrix}
  1_{n}\otimes \alpha'  &          0      \\
    0               &   1_{n}\otimes \alpha'
\end{bmatrix},
\begin{bmatrix}
     1_{n}\otimes \beta'   &    0      \\
    0                 &   1_{n}\otimes \beta'
\end{bmatrix}
)$$
\end{definition}

We now recall the definition of the multiplicative tensor product of two maps.\\
Let $X_{f}=(\phi,\psi)$, $X_{f}'=(\phi',\psi')$ and $X_{f}"=(\phi",\psi")$ be objects of $MF(K[[x]],f)$ respectively of sizes $n, n'$ and $n"$. Let $X_{g}=(\sigma,\rho)$, $X_{g}'=(\sigma',\rho')$ and $X_{g}"=(\sigma",\rho")$ be objects of $MF(K[[y]],g)$ respectively of sizes $m, m'$ and $m"$.

\begin{definition}\label{defn zeta is a bifuntor} \cite{fomatati2019multiplicative}
For morphisms $\zeta_{f}=(\alpha_{f}, \beta_{f}): X_{f}=(\phi,\psi) \rightarrow X_{f}'=(\phi',\psi')$  and $\zeta_{g}=(\alpha_{g}, \beta_{g}): X_{g}=(\sigma,\rho) \rightarrow X_{g}'=(\sigma',\rho')$ respectively in $MF(K[[x]],f)$ and $MF(K[[y]],g)$, we define $\zeta_{f}\widetilde{\otimes} \zeta_{g}: X_{f}\widetilde{\otimes}X_{g} =(\phi,\psi)\widetilde{\otimes} (\sigma,\rho)\rightarrow X_{f}'\widetilde{\otimes}X_{g}' =(\phi',\psi')\widetilde{\otimes} (\sigma',\rho')$
by
$$(
\begin{bmatrix}
  \alpha_{f}\otimes \alpha_{g}  &          0      \\
    0               &   \alpha_{f}\otimes \alpha_{g}
\end{bmatrix},
\begin{bmatrix}
     \beta_{f}\otimes \beta_{g}  &    0      \\
    0                 &   \beta_{f}\otimes \beta_{g}
\end{bmatrix}
)$$

\end{definition}
\begin{theorem} \label{thm: mult tens prdt is a bifunctor} \cite{fomatati2019multiplicative}
  $\widetilde{\otimes}$ is a bifunctor.
\end{theorem}
A variant of $\widetilde{\otimes}$ was found in \cite{fomatati2021tensor}.
%
%
%
%
\section{A comparison of $\widehat{\otimes}$ and $\widetilde{\otimes}$, and a study of the category $(MF(1),\widetilde{\otimes})$}
In this section, we compare $\widehat{\otimes}$ and $\widetilde{\otimes}$ and study the category $(MF(1),\widetilde{\otimes})$. \\
The Syzygy property (cf. subsection \ref{subsec: comparison of the two tensor prdts} below) will help to find some differences between these two operations.
 Moreover, we will observe that the multiplicative tensor product of two objects of $MF(1)$ is still an object of $MF(1)$ whereas the tensor product $\widehat{\otimes}$ of any two matrix factorizations of a power series $f$ is not a matrix factorization of $f$ (not even for $f=1$). This will motivate the study of $(MF(1),\widetilde{\otimes})$. Is it a monoidal category? or a generalization of this notion? \\
 We will define the concept of \textit{one-step connected category} and prove that there is a one-step connected subcategory $(\mathcal{T},\widetilde{\otimes})$ of $(MF(1),\widetilde{\otimes})$ which is a semi-unital semi-monoidal category. This is particularly interesting because the concept of semi-unital semi-monoidal category was recently conceived in \cite{abuhlail2013semiunital} and the example provided in that paper (cf. theorem 5.11 of \cite{abuhlail2013semiunital}) required a considerable amount of set-up. But in this section, the example (cf. theorem \ref{Thm: T is a semi-unital semi-monoidal catg}) we give requires a smaller amount of set-up. \\
 Furthermore, we will define the concept of \textit{right pseudo-monoidal} category and prove that the category $(MF(1),\widetilde{\otimes})$ is an example of this concept.\\
First, recall \cite{fomatati2019multiplicative} that
  a $(1,0)$-matrix is a matrix whose entries belong to the set $\{0,1\}$.
We chose the terminology $(1,0)-$matrix instead of $(0,1)-$matrix because some authors use the terminology $(0,1)-$matrix to refer to what we call here $(1,0)-$matrix with some additional
conditions.

\begin{definition}
  A category is said to be a \textbf{one-step connected category} if for every two objects of the category, there exists a nonzero morphism between them.
\end{definition}

\subsection{A comparison of $\widehat{\otimes}$ and $\widetilde{\otimes}$} \label{subsec: comparison of the two tensor prdts}
$\widehat{\otimes}$ and $\widetilde{\otimes}$ are different at several levels. First of all,
the Syzygy $\Omega$ helps in pointing out some differences between these two operations.\\
First recall the definition of the Syzygy\footnote{We use this word because that is the name (cf. \cite{yoshino1998tensor}) given to the operator $\Omega$ we are going to use in Prop. \ref{prop: Syzygy}.} $\Omega((\phi,\psi)):=(\psi,\phi)$ where $(\psi,\phi)$ is a matrix factorization of a power series $f$.\\
We now want to state a Syzygy property for $\widetilde{\otimes}$. In \cite{yoshino1998tensor}, a Syzygy property was proved for $\widehat{\otimes}$, the \textit{tensor product of matrix factorization} (cf. subsection \ref{subsec: Yoshino tens prod and variant}). It was proved that $X\widehat{\otimes}X'=\Omega(X)\widehat{\otimes}\Omega(X')$ and $X\widehat{\otimes}\Omega(X')\cong \Omega(X\widehat{\otimes}X') \cong \Omega(X)\widehat{\otimes}X'$. But the Syzygy property that holds for $\widetilde{\otimes}$ is totally different. It shows that the functor $\Omega$ is "linear" with respect to the operation $\widetilde{\otimes}$.
\begin{proposition} \cite{fomatati2019multiplicative} (Syzygy property) \label{prop: Syzygy}\\
There is an identity $\Omega(X\widetilde{\otimes}X')=\Omega(X)\widetilde{\otimes}\Omega(X')$.
\end{proposition}

For instance, observe that in general, unlike with $\widehat{\otimes}$; $\Omega(X)\widetilde{\otimes}X'\ncong X\widetilde{\otimes}\Omega(X')$ as can be easily checked computationally. \\
%
It is easy to verify that in general, $\Omega(X)\widetilde{\otimes}X'\ncong X\widetilde{\otimes}\Omega(X')$. It is also easy to verify that unlike with $\widehat{\otimes}$, $X \widetilde{\otimes}X'\neq \Omega(X)\widetilde{\otimes}\Omega(X')$.\\
Moreover, from the definitions of $\widehat{\otimes}$ (cf. subsection \ref{subsec: Yoshino tens prod and variant}) and $\widetilde{\otimes}$ (cf. definition \ref{defn of the multiplicative tensor product}), we immediately see some similarities and differences. For example, given two matrix factorizations $X_{f}$ of $f\in K[[x]]$ of size $n$ and $X_{g}$ of $g\in K[[y]]$ of size $m$, though $X_{f}\widehat{\otimes}X_{g}$ and $X_{f}\widetilde{\otimes} X_{g}$ are both of size  $2nm$, they are objects of two different categories namely $MF(f+g)$ and $MF(fg)$. Even if we consider two objects of the same category, say $X_{f}$ and $Y_{f}$ of a nonzero power series $f\in K[[x]]$, $X_{f}\widehat{\otimes}Y_{f}$ (respectively $X_{f}\widetilde{\otimes} Y_{f}$) will be an object not in $MF(f)$ but instead of a different category namely $MF(f+f)$ (respectively $MF(f\cdot f)=MF(f^{2})$). Now, there is a striking difference between the two tensor products when $f=1$. In fact, if $f=1$, then $X_{f}\widehat{\otimes}Y_{f}\in MF(1+1)\neq MF(1)$ but $X_{f}\widetilde{\otimes} Y_{f}\in MF(1\cdot 1)=MF(1)$. That is, the multiplicative tensor product of two objects of $MF(1)$ is still an object of $MF(1)$. This motivates the study of $(MF(1),\widetilde{\otimes})$ to know whether it is a monoidal category or a generalization of this notion.

\subsection{An application of $\widetilde{\otimes}$: A semi-unital semi-monoidal subcategory of $MF(1)$} \label{subsec: a semi-unital semi-monoidal subcatg of MF(1)}
We prove that $MF(1)$ has a one-step connected subcategory which is a semi-unital semi-monoidal category. \\
Objects of $MF(1)$ are of the form $(M_{n}, M_{n}^{-1})$ where $M_{n}$ is an $n\times n$ matrix, $n\in \mathbb{N}-\{0\}$. Morphisms are pairs of matrices such that a certain diagram commutes (cf. subsection \ref{subsec: Yoshino tens prod and variant}).

\begin{theorem} \label{Thm: T is a semi-unital semi-monoidal catg}
There is a one-step connected subcategory of $MF(1)$ which is a semi-unital semi-monoidal category.
\end{theorem}

\begin{proof}
We extract a one-step connected subcategory of $MF(1)$ which is a semi-unital semi-monoidal category. We will call it $\mathcal{T}$.\\
$\bullet$ Objects of $\mathcal{T}$ are of the form $e^{n}$ where $n\in \mathbb{N}-\{0\}$. We characterize these objects. \\
$e^{1}=e=(1,1)$, $e^{2}=e\widetilde{\otimes}e=\left(\left(\begin{matrix}
1\otimes 1 &  0\\
0  &  1\otimes 1\\
\end{matrix}\right), \left(\begin{matrix}
1\otimes 1 &  0\\
0  &  1\otimes 1\\
\end{matrix}\right)\right)=(I_{2},I_{2})$
$$e^{3}=e\widetilde{\otimes} e^{2}= (1,1)\widetilde{\otimes} \left(\left(\begin{matrix}
1\otimes 1 &  0\\
0  &  1\otimes 1\\
\end{matrix}\right), \left(\begin{matrix}
1\otimes 1 &  0\\
0  &  1\otimes 1\\
\end{matrix}\right)\right)$$
$$=\left(\begin{matrix}1\otimes \left(\begin{matrix}
1\otimes 1 &  0\\
0  &  1\otimes 1\\
\end{matrix}\right) & 0\\
0   &  1\otimes \left(\begin{matrix}
1\otimes 1 &  0\\
0  &  1\otimes 1\\
\end{matrix}\right)\end{matrix}\right),\left(\begin{matrix}1\otimes \left(\begin{matrix}
1\otimes 1 &  0\\
0  &  1\otimes 1\\
\end{matrix}\right) & 0\\
0   &  1\otimes \left(\begin{matrix}
1\otimes 1 &  0\\
0  &  1\otimes 1\\
\end{matrix}\right)\end{matrix}\right)=(I_{4},I_{4})$$
It is easy to see that in general, $e^{n}=(I_{2^{n-1}}, I_{2^{n-1}})$.\\
We give a notation before defining morphisms between two objects of $\mathcal{T}$.\\
\begin{notation}
  We will denote by $0_{n,m}$ the zero matrix of size $n\times m$ whenever there would be a risk of confusion on the size of the zero matrix in the context under consideration. Otherwise, we will simply write $0$. 
\end{notation}
$\bullet$ Morphisms of $\mathcal{T}$ are defined as in $MF(1)$ (cf. subsection \ref{subsec: Yoshino tens prod and variant}), but with some restrictions.
We now define what a morphism is in $\mathcal{T}$.  \\
\textbf{Discussion $\sharp $}:\\
First recall that a permutation matrix is a square matrix obtained from the same size identity matrix by a permutation of rows.\\
 For $m,p\in \mathbb{N}-\{0\}$, recall that a morphism $e^{m}\rightarrow e^{p}$ in $MF(1)$  is a pair of matrices $(\delta, \beta)$ such that the following diagram commutes:
 $$\xymatrix@ R=0.6in @ C=.75in{K[[x]]^{2^{m-1}} \ar[r]^{I_{2^{m-1}}} \ar[d]_{\delta} &
K[[x]]^{2^{m-1}} \ar[d]^{\beta} \ar[r]^{I_{2^{m-1}}} & K[[x]]^{2^{m-1}}\ar[d]^{\delta \;\;\;\;\;\;\;\;\;\;(\bigstar")}\\
K[[x]]^{2^{p-1}} \ar[r]^{I_{2^{p-1}}} & K[[x]]^{2^{p-1}}\ar[r]^{I_{2^{p-1}}} & K[[x]]^{2^{p-1}}}$$
That is,
$$ (S) \begin{cases}
 \delta I_{2^{m-1}}=I_{2^{p-1}}\beta  \\
 I_{2^{p-1}}\delta= \beta I_{2^{m-1}}
\end{cases}$$
It follows from $(S)$, that a morphism $e^{m}\rightarrow e^{p}$ in $MF(1)$  is a pair of matrices $(\delta, \beta)$ with $\delta=\beta$. This does not impose any restrictions on the entries of $\delta$ or $\beta$, the entries could be anything provided we have the equality $\delta=\beta$.\\
\\ But, we define a morphism $e^{m}\rightarrow e^{p}$ in $\mathcal{T}$ to be a pair of matrices $(\delta, \beta)$ such that $\delta=\beta$ is a $(1,0)-$matrix of size $2^{p-1}\times 2^{m-1}$ \textit{with at most one nonzero entry in each row and each column}. This restriction will ensure that the composition of two morphisms in $\mathcal{T}$ is again a morphism in $\mathcal{T}$.\\
Thus, for example we could have the following values of $\delta$ and $\beta$ for the pair $(\delta, \beta)$ to be a morphism in $\mathcal{T}$:\\

  $\delta=\beta=\begin{cases}
(I_{2^{p-1}},0_{2^{p-1},2^{m-1}-2^{p-1}}), \,if\, m>p\\

\begin{pmatrix}
          I_{2^{m-1}}  \\
0_{2^{p-1}-2^{m-1},2^{m-1}}
\end{pmatrix},\, if \,m<p  \\
Z\, if\,m=p \,where\,Z \,is\,a\,2^{m-1}\times 2^{m-1}\, permutation\, matrix
\end{cases}$
\\
 From the twin square diagram $\bigstar"$, it is clear that $\delta$ and $\beta$ should both be of size $2^{p-1}\times 2^{m-1}$.
The fact that we actually have a morphism from $e^{m}$ to $e^{p}$ for the above values of $\delta=\beta$ is obvious from diagram $\bigstar"$.\\
Discussion $\sharp$ actually gives us a sufficient condition on a pair $(\delta,\beta)$ to be a nonzero morphism in $\mathcal{T}$.\\

It is not difficult to see that $\mathcal{T}$ is a subcategory of $MF(1)$.\\
In fact, for every pair of morphisms $\zeta$ and $\zeta'$ in $hom(\mathcal{T})$, the composite $\zeta \circ \zeta'$ is in $hom(\mathcal{T})$ whenever it is defined. In fact, $\zeta$ and $\zeta'$ by definition are pairs of $(1,0)-$matrices such that in each matrix; each column and each row has at most one nonzero entry. It then follows that the composition of such matrices will yield another $(1,0)-$matrix in which each column and each row would have at most one nonzero entry, whence we will still have a morphism of $\mathcal{T}$. \\

Moreover, $\mathcal{T}$ is a one-step connected category because between any two objects of $\mathcal{T}$ say $e^{m}$ and $e^{p}$, there exists a nonzero morphism as can be seen from discussion $\sharp$.\\

We now proceed to prove that $\mathcal{T}$ is a semi-unital semi-monoidal category.\\
$\bullet$ We first prove that $(\mathcal{T}, \widetilde{\otimes})$ is a semi-monoidal category (cf. definition \ref{defn semi-monoidal category}):\\
- The fact that $\widetilde{\otimes}: MF(1)\times MF(1) \rightarrow MF(1)$ is a bifunctor follows from theorem 3.1 \cite{fomatati2019multiplicative} by replacing $f$ and $g$ by the constant power series $1$ and by letting $x=y$.\\
- There is a natural isomorphism $\alpha$ from the functor $((-)\widetilde{\otimes}(-))\widetilde{\otimes} (-): MF(1)\times MF(1)\times MF(1) \rightarrow (MF(1)\times MF(1))\times MF(1)$ to the functor $(-)\widetilde{\otimes}((-)\widetilde{\otimes} (-)): MF(1)\times MF(1)\times MF(1) \rightarrow MF(1)\times (MF(1)\times MF(1))$ with components $\alpha_{a,b,c}: (a\widetilde{\otimes} b)\widetilde{\otimes}c \rightarrow a\widetilde{\otimes} (b\widetilde{\otimes}c)$, where $a,b\,and\,c$ are matrix factorizations of $1$ in $\mathcal{T}$.\\
Let $a,b\,c,a',b',and\,c'$ be objects of $\mathcal{T}$. Let $f: a \rightarrow a'$, $g: b \rightarrow b'$ and $h: c \rightarrow c'$ be maps in $\mathcal{T}$. We show that the following diagram commutes:\\
  \[
\xymatrix{
(a\widetilde{\otimes} b)\widetilde{\otimes}c\ar[r]^{\alpha_{a,b,c}} \ar[d]_{(f\widetilde{\otimes}g) \widetilde{\otimes} h} & a\widetilde{\otimes} (b\widetilde{\otimes}c) \ar[d]_{f\widetilde{\otimes}(g \widetilde{\otimes} h)}\\
(a'\widetilde{\otimes} b')\widetilde{\otimes}c'\ar[r]_{\alpha_{a',b',c'}} & a'\widetilde{\otimes} (b'\widetilde{\otimes}c')}
\]
  \\
  i.e., $f\widetilde{\otimes}(g \widetilde{\otimes} h) \circ \alpha_{a,b,c}=\alpha_{a',b',c'} \circ (f\widetilde{\otimes}g) \widetilde{\otimes} h$ $\cdots (E')$ \\
  In fact, the matrices representing $\alpha_{a,b,c}$ and $\alpha_{a',b',c'}$ are identity matrices. Besides, the tensor product of maps is associative. Thus, $(E')$ holds. That is $\alpha$ is a natural transformation. Moreover, for all $a,b$ and $c$; $\alpha_{a,b,c}$ is an equality and so, it is an isomorphism. Hence, $\alpha$ is a natural isomorphism. \\
Next, let us show that the pentagonal diagram of definition \ref{defn semi-monoidal category} commutes for all $a,b,c,d\,\in \,MF(1)$
$\xymatrix @ R=0.4in @ C=1.0in{
a\widetilde{\otimes} (b\widetilde{\otimes} (c\widetilde{\otimes} d)) \ar[rr]^{1_{a}\widetilde{\otimes} \alpha}\ar[d]_\alpha &&
a\widetilde{\otimes} ((b\widetilde{\otimes} c)\widetilde{\otimes} d) \ar[d]^\alpha &\\
(a\widetilde{\otimes} b)\widetilde{\otimes} (c\widetilde{\otimes} d) \ar[rd]_\alpha &&
(a\widetilde{\otimes}(b\widetilde{\otimes} c))\widetilde{\otimes} d \ar[ld]^{\alpha \widetilde{\otimes} 1_{d}} &\\
&((a\widetilde{\otimes} b)\widetilde{\otimes} c)\widetilde{\otimes} d} $

Since all the maps linking the vertices of the pentagon are identity maps, this diagram must commute. In fact, we know that $\alpha$ (the associator) is an identity map. Moreover, since the pair of matrices making up $\alpha$ are identity matrices, it follows from the definition \ref{defn zeta is a bifuntor} (of the multiplicative tensor product of two maps) that $1_{a}\widetilde{\otimes} \alpha$ and $\alpha \widetilde{\otimes} 1_{d}$ are also identity maps.

  Therefore, $\mathcal{T}$ is a semi-monoidal category.\\\\
$\bullet$ Next, we prove that $\mathcal{T}$ is a semi-unital semi-monoidal category. To that end, we need to find a semi-unit in the semi-monoidal category $\mathcal{T}$.\\
\textbf{Claim:} $e=(1,1)$ is a semi-unit in $\mathcal{T}$.\\
From the definition of a semi-unit (cf. definition \ref{defn semi-unital semimonoidal catg}), we need to find a natural transformation $\gamma : (-)=G \rightarrow F=e\widetilde{\otimes} (-)$, where G is the identity endofunctor on $\mathcal{T}$ and $F$ is an endofunctor on $\mathcal{T}$, such that $F(a)=e\widetilde{\otimes} a$. Components of $\gamma$ are:\\
$\gamma_{a}: a \rightarrow e\widetilde{\otimes} a$, where $a=e^{p}=(\phi,\psi)$ is an object of $\mathcal{T}$ of size $n_{1}$, that is $\phi=I_{2^{p-1}}=\psi$ with $n_{1}=2^{p-1}$. We have:  $e\widetilde{\otimes} a=\left(\begin{bmatrix}
I_{2^{p-1}}  &    0      \\
    0   &  I_{2^{p-1}}
\end{bmatrix}, \begin{bmatrix}
I_{2^{p-1}}  &    0      \\
    0   &  I_{2^{p-1}}
\end{bmatrix}\right)=(I_{2^{p}}, I_{2^{p}})$ is of size $2^{p}=2n_{1}$. \\
The family of morphisms $\gamma$ should satisfy the following two requirements:
\begin{enumerate}
  \item For each $a\in\,Ob(\mathcal{T})$, $\gamma_{a}$ should be a morphism in $\mathcal{T}$.\\
 Since $a$ and $e\widetilde{\otimes} a$ are objects of $\mathcal{T}$ which is a one-step connected category, we let $\gamma_{a}$ be the nonzero morphism between $a$ and $e\widetilde{\otimes} a$ such that:  $\gamma_{a}=(\delta',\beta')=(\delta',\delta')= ((I_{n_{1}},0)^{t}, (I_{n_{1}},0)^{t})$, where $t$ is the operation of taking the transpose, $0$ is the zero $n_{1}\times n_{1}$ matrix. $\gamma_{a}$ is clearly a morphism in $\mathcal{T}$ as discussed under discussion $\sharp$.

\item Naturality of $\gamma$:\\
 Let $b=(\phi',\psi')$ be a matrix factorization in $\mathcal{T}$ of size $n_{2}$ and let $\mu= (\alpha_{\mu}, \beta_{\mu}): a\rightarrow b$ be a map of matrix factorizations. It is easy
 to see that $\alpha_{\mu}$ and $\beta_{\mu}$ are each of size $n_{2}\times n_{1}$.
  The following diagram should commute:\\
  \[
\xymatrix{
 a\ar[r]^{\gamma_{a}} \ar[d]_{\mu} & e\widetilde{\otimes} a \ar[d]_{e\widetilde{\otimes}\mu}\\
 b\ar[r]_{\gamma_{b}} & e\widetilde{\otimes} b}
\]
  \\
  i.e., $e\widetilde{\otimes}\mu \circ \gamma_{a}=\gamma_{b} \circ \mu$ $\cdots (E'')$\\
  We know that $e\widetilde{\otimes} a$ is of size $2n_{1}$ since $a$ is of size $n_{1}$. We also know that $\gamma_{b}=[(I_{n_{2}},0)^{t},(I_{n_{2}},0)^{t}]$. Now, by definition of composition of two morphisms in $\mathcal{T}$, the right hand side of equality $(E'')$ becomes:\\
  $\gamma_{b}\circ \mu= [(I_{n_{2}},0)^{t},(I_{n_{2}},0)^{t^{}}]\circ (\alpha_{\mu}, \beta_{\mu})=[(I_{n_{2}},0)^{t}\alpha_{\mu}, (I_{n_{2}},0)^{t}\beta_{\mu}]=[(\alpha_{\mu},0)^{t},(\beta_{\mu},0)^{t}]$$\cdots \natural'$\\
  $0$ in $[(\alpha_{\mu},0)^{t},(\beta_{\mu},0)^{t}]$ is the $n_{2}\times n_{1}$ zero matrix. \\
  As for the left hand side of $(E'')$, first recall that $\gamma_{a}=[(I_{n_{1}},0)^{t},(I_{n_{1}},0)^{t}]$, (where $0$ is the zero $n_{1}\times n_{1}$ matrix) and by definition \ref{defn zeta is a funtor (right)} of the multiplicative tensor product, we know that $e\widetilde{\otimes} \mu=(1,1)\widetilde{\otimes}(\alpha_{\mu}, \beta_{\mu})=(\begin{bmatrix}
1 \otimes \alpha_{\mu}  &    0      \\
    0                 &  1 \otimes \alpha_{\mu}
\end{bmatrix},
\begin{bmatrix}
1 \otimes \beta_{\mu}  &    0      \\
    0                 &  1 \otimes \beta_{\mu}
\end{bmatrix})=(\begin{bmatrix}
\alpha_{\mu} &    0      \\
    0        &  \alpha_{\mu}
\end{bmatrix},
\begin{bmatrix}
\beta_{\mu}  &    0      \\
    0          &  \beta_{\mu}
\end{bmatrix})$ \\
So, $e\widetilde{\otimes} \mu \circ \gamma_{a} = (\begin{bmatrix}
\alpha_{\mu} &    0      \\
    0        &  \alpha_{\mu}
\end{bmatrix},
\begin{bmatrix}
\beta_{\mu}  &    0      \\
    0          &  \beta_{\mu}
\end{bmatrix})\circ [(I_{n_{1}},0)^{t},(I_{n_{1}},0)^{t}]\\
=(\begin{bmatrix}
\alpha_{\mu} &    0      \\
    0        &  \alpha_{\mu}
\end{bmatrix}(I_{n_{1}},0)^{t},\begin{bmatrix}
\beta_{\mu}  &    0      \\
    0          &  \beta_{\mu}
\end{bmatrix}(I_{n_{1}},0)^{t})=[(\alpha_{\mu},0)^{t},(\beta_{\mu},0)^{t}] \cdots \natural\natural'$.\\
From $\natural'$ and $\natural\natural'$, we see that equality $(E'')$ holds.
\end{enumerate}

Hence $\gamma$ is a natural transformation.\\\\
The next step towards proving that $\mathcal{T}$ is a semi-unital semi-monoidal category is to prove
that there is an isomorphism of functors $e \widetilde{\otimes} (-) \cong (-) \widetilde{\otimes} e$, i.e., there is a natural isomorphism $l_{a}:e \widetilde{\otimes} (a) \cong (a) \widetilde{\otimes} e$ with inverse $q_{a}$, for each object $a$ of $\mathcal{T}$ such that $l_{e}=q_{e}$ and the following diagrams are commutative for all objects $a$ and $b$ of $\mathcal{T}$.\\

  $$\xymatrix@ R=0.6in @ C=.75in{(e\widetilde{\otimes} a)\widetilde{\otimes} b\ar[r]^{\alpha_{e,a,b}} \ar[d]_{l_{a}\widetilde{\otimes} b} &
e\widetilde{\otimes} (a\widetilde{\otimes} b) \ar[r]^{l_{a\widetilde{\otimes} b}} & (a\widetilde{\otimes} b)\widetilde{\otimes} e\ar[d]^{\alpha_{a,b,e}\,\,\,\,\,\,\,\,\,\,\,\,\,\,\,\,(1)}\\
(a\widetilde{\otimes} e)\widetilde{\otimes} b \ar[r]^{\alpha_{a,e,b}} & a\widetilde{\otimes} (e\widetilde{\otimes} b)\ar[r]^{a \widetilde{\otimes} l_{b}} & a\widetilde{\otimes} (b\widetilde{\otimes} e)}$$

\[
\xymatrix{
a\widetilde{\otimes} b\ar[rr]^{\gamma_{a}\widetilde{\otimes} b} \ar[rd]_{\gamma_{a\widetilde{\otimes} b}} &&
(e\widetilde{\otimes} a) \widetilde{\otimes} b \ar[ld]^{\cong \,\,\,\,\,\,\,\,\,\,\,\,\,\,\,\,(2)}\\
&e\widetilde{\otimes} (a \widetilde{\otimes} b)}
\]
\[
\xymatrix{
a\widetilde{\otimes} b\ar[rr]^{a \widetilde{\otimes}\gamma_{b}} \ar[rd]_{\gamma_{a\widetilde{\otimes} b}} &&
a\widetilde{\otimes} (e \widetilde{\otimes} b) \ar[ld]^{\cong \,\,\,\,\,\,\,\,\,\,\,\,\,\,\,\,(3)}\\
&e\widetilde{\otimes} (a \widetilde{\otimes} b)}
\]
Before we define $l_{a}$, observe that  $e \widetilde{\otimes} (a)= (a) \widetilde{\otimes} e$.\\
$\clubsuit$ We define $l_{a}: e \widetilde{\otimes} (a) \rightarrow (a) \widetilde{\otimes} e$ to be the pair of matrices $(I_{2n_{1}}, I_{2n_{1}})=(I_{2^{p}}, I_{2^{p}})$ where $a=e^{p}$ is of size $n_{1}$. From discussion $\sharp$, it follows that $l_{a}$ is a morphism in $\mathcal{T}$.\\
$\clubsuit$ Naturality of $l$:\\

Let $b=(\phi',\psi')$ be a matrix factorization of size $n_{2}$ and let $\mu'= (\alpha_{\mu'}, \beta_{\mu'}): a\rightarrow b$ be a map of matrix factorizations. It is easy\footnote{By drawing the twin diagram that has to commute with $(\alpha_{\mu'}, \beta_{\mu'})$, we see the sizes of $\alpha_{\mu'}$ and $\beta_{\mu'}$. } to see that $\alpha_{\mu'}$ and $\beta_{\mu'}$ are each of size $n_{2}\times n_{1}$.
  The following diagram should commute:\\
  \[
\xymatrix{
 e\widetilde{\otimes} a\ar[r]^{l_{a}} \ar[d]_{e \widetilde{\otimes} \mu'} & a\widetilde{\otimes} e \ar[d]_{\mu' \widetilde{\otimes}e}\\
 e\widetilde{\otimes} b\ar[r]_{l_{b}} & b\widetilde{\otimes} e}
\]
  \\
  i.e., $\mu' \widetilde{\otimes} e \circ l_{a}=l_{b} \circ e\widetilde{\otimes} \mu'$ $\cdots (E''')$\\
  Since $l_{a}$ and $l_{b}$ are just pairs of identity matrices, it suffices to show that
  $\mu' \widetilde{\otimes} e = e\widetilde{\otimes} \mu'$.  \\
By definition \ref{defn zeta is a funtor (right)} of the multiplicative tensor product, we know that $e\widetilde{\otimes} \mu'=(1,1)\widetilde{\otimes}(\alpha_{\mu'}, \beta_{\mu'})=(\begin{bmatrix}
1 \otimes \alpha_{\mu'}  &    0      \\
    0                 &  1 \otimes \alpha_{\mu'}
\end{bmatrix},
\begin{bmatrix}
1 \otimes \beta_{\mu'}  &    0      \\
    0                 &  1 \otimes \beta_{\mu'}
\end{bmatrix})=(\begin{bmatrix}
\alpha_{\mu'} &    0      \\
    0        &  \alpha_{\mu'}
\end{bmatrix},
\begin{bmatrix}
\beta_{\mu'}  &    0      \\
    0          &  \beta_{\mu'}
\end{bmatrix})$ \\
And we also have by definition \ref{defn zeta tensor tilda X'} of the multiplicative tensor product, that \\ $\mu' \widetilde{\otimes} e=(\alpha_{\mu'}, \beta_{\mu'})\widetilde{\otimes}(1,1)=(\begin{bmatrix}
\alpha_{\mu'} \otimes 1  &    0      \\
    0                 &  \alpha_{\mu'} \otimes 1
\end{bmatrix},
\begin{bmatrix}
 \beta_{\mu'}\otimes 1  &    0      \\
    0                 &  \beta_{\mu'} \otimes 1
\end{bmatrix})=(\begin{bmatrix}
\alpha_{\mu'} &    0      \\
    0        &  \alpha_{\mu'}
\end{bmatrix},
\begin{bmatrix}
\beta_{\mu'}  &    0      \\
    0          &  \beta_{\mu'}
\end{bmatrix}).$ \\
  Thus $e\widetilde{\otimes} \mu'=\mu' \widetilde{\otimes} e $, so $l$ is a natural transformation.\\
$\clubsuit$ $l_{a}:e \widetilde{\otimes} (a) \cong (a) \widetilde{\otimes} e$ is a natural isomorphism. In fact, $l: e \widetilde{\otimes} (a) =(a) \widetilde{\otimes} e$ and its inverse $q:(a) \widetilde{\otimes} e= e \widetilde{\otimes} (a)$ is clearly such that $l_{e}=q_{e}$.\\

\textbf{Commutativity of diagrams $(1), (2) \,and\,(3)$:}\\
For diagram $(1)$, it would commute if $\alpha_{a,b,e} \circ l_{a\widetilde{\otimes}b} \circ \alpha_{e,a,b}=a\widetilde{\otimes} l_{b} \circ \alpha_{a,e,b} \circ l_{a} \widetilde{\otimes} b \,\,\,\,\,\,\,\,\,\,\,\, \cdots (\flat)$.\\
We show that all the maps involved in equality $(\flat)$ are identities. For all objects $x,y,z \,in\, Ob(\mathcal{T})$, we clearly have by the definitions of $\alpha_{x,y,z}$ and $l_{x}$ that they are identity maps. We now show that the other maps involved in diagram $(1)$ are identity maps.\\
Since $a$ is of size $n_{1}$, we have $l_{a}=(I_{2n_{1}}, I_{2n_{1}})$. Let $b=(\phi',\psi')$ be of size $n_{2}$. By definition \ref{defn zeta tensor tilda X'}, $l_{a}\widetilde{\otimes} b =(I_{2n_{1}}, I_{2n_{1}})\widetilde{\otimes}(\phi',\psi')=(\begin{bmatrix}
I_{2n_{1}} \otimes I_{n_{2}}  &    0      \\
    0                 &  I_{2n_{1}} \otimes I_{n_{2}}
\end{bmatrix},
\begin{bmatrix}
 I_{2n_{1}} \otimes I_{n_{2}}  &    0      \\
    0                 &  I_{2n_{1}} \otimes I_{n_{2}}
\end{bmatrix})=(\begin{bmatrix}
I_{2n_{1}n_{2}} &    0      \\
    0        &  I_{2n_{1}n_{2}}
\end{bmatrix},
\begin{bmatrix}
I_{2n_{1}n_{2}} &    0      \\
    0        &  I_{2n_{1}n_{2}}
\end{bmatrix})=(I_{4n_{1}n_{2}}, I_{4n_{1}n_{2}})  \cdots \ddag$\\
Hence $l_{a}\widetilde{\otimes} b$ is an identity map as expected.\\

Similarly using definition \ref{defn zeta is a funtor (right)}, we prove that $a\widetilde{\otimes} l_{b}$ is an identity map. Let $a=(\phi,\psi)$ be of size $n_{1}$ and $b$ be as above. Then, $a\widetilde{\otimes} l_{b} =(\phi,\psi)\widetilde{\otimes}(I_{2n_{2}}, I_{2n_{2}})\\
=(\begin{bmatrix}
I_{n_{1}} \otimes I_{2n_{2}}  &    0      \\
    0                 &  I_{n_{1}} \otimes I_{2n_{2}}
\end{bmatrix},
\begin{bmatrix}
I_{n_{1}} \otimes I_{2n_{2}}  &    0      \\
    0                 &  I_{n_{1}} \otimes I_{2n_{2}}
\end{bmatrix})=(\begin{bmatrix}
I_{2n_{1}n_{2}} &    0      \\
    0        &  I_{2n_{1}n_{2}}
\end{bmatrix},
\begin{bmatrix}
I_{2n_{1}n_{2}} &    0      \\
    0        &  I_{2n_{1}n_{2}}
\end{bmatrix})=(I_{4n_{1}n_{2}}, I_{4n_{1}n_{2}})    \cdots \ddag'$ \\

$\ddag$ and $\ddag'$ show that $a\widetilde{\otimes} l_{b}=l_{a}\widetilde{\otimes} b$.\\
It is easy to see that all the other maps involved in diagram $(1)$ are equal to $(I_{4n_{1}n_{2}}, I_{4n_{1}n_{2}})$.
\\
So diagram $(1)$ is commutative.\\\\
Next, we show that diagram $(2)$ commutes. To this end, we need to find an isomorphism $\zeta: (e\widetilde{\otimes} a) \widetilde{\otimes} b \rightarrow e\widetilde{\otimes} (a \widetilde{\otimes} b)$,
such that $\zeta \circ \gamma_{a}\widetilde{\otimes} b = \gamma_{a\widetilde{\otimes} b}$.\\
Now, we know that $\gamma_{a}=[(I_{n_{1}}, 0)^{t}, (I_{n_{1}}, 0)^{t}]$, where $0$ is the $n_{1}\times n_{1}$ zero matrix and $b=(\phi',\psi')$ is of size $n_{2}$. Hence by definition \ref{defn zeta tensor tilda X'} $\gamma_{a}\widetilde{\otimes} b= [(I_{n_{1}}, 0)^{t}, (I_{n_{1}}, 0)^{t}]\widetilde{\otimes}(\phi',\psi')=(\begin{bmatrix}
(I_{n_{1}}, 0_{n_{1},n_{1}})^{t} \otimes I_{n_{2}}  &    0_{2n_{1}n_{2},n_{1}n_{2}}      \\
    0_{2n_{1}n_{2},n_{1}n_{2}}                 &  (I_{n_{1}}, 0_{n_{1},n_{1}})^{t} \otimes I_{n_{2}}
\end{bmatrix},
\begin{bmatrix}
 (I_{n_{1}}, 0_{n_{1},n_{1}})^{t} \otimes I_{n_{2}}  &    0_{2n_{1}n_{2},n_{1}n_{2}}      \\
    0_{2n_{1}n_{2},n_{1}n_{2}}                 &  (I_{n_{1}}, 0_{n_{1},n_{1}})^{t} \otimes I_{n_{2}}
\end{bmatrix})=\\
(\begin{bmatrix}
(I_{n_{1}n_{2}}, 0_{n_{1}n_{2}})^{t} &    0_{2n_{1}n_{2},n_{1}n_{2}}      \\
    0_{2n_{1}n_{2},n_{1}n_{2}}        &  (I_{n_{1}n_{2}}, 0_{n_{1}n_{2}})^{t}
\end{bmatrix},
\begin{bmatrix}
(I_{n_{1}n_{2}}, 0_{n_{1}n_{2}})^{t} &    0_{2n_{1}n_{2},n_{1}n_{2}}      \\
    0_{2n_{1}n_{2},n_{1}n_{2}}        &  (I_{n_{1}n_{2}}, 0_{n_{1}n_{2}})^{t}
\end{bmatrix}).$ $\cdots \dag$\\
Next, since $a\widetilde{\otimes} b$ is an object of size $2n_{1}n_{2}$, we obtain from the way $\gamma$ is defined that
$\gamma_{a\widetilde{\otimes} b}= ((I_{2n_{1}n_{2}}, 0_{2n_{1}n_{2}})^{t}, (I_{2n_{1}n_{2}}, 0_{2n_{1}n_{2}})^{t}$) $\cdots \dag'$ \\
From $\dag$ and $\dag'$, we see that $\gamma_{a\widetilde{\otimes} b}$ and $\gamma_{a}\widetilde{\otimes} b$ are both $(1,0)-$matrices with the same number of rows and columns. Moreover, they have the same number of $1$s and each of these $1$s is the only nonzero entry in its row and in its column. Simply put, the matrix from $\gamma_{a\widetilde{\otimes} b}$ is (row-)permutation equivalent to the matrix in $\gamma_{a}\widetilde{\otimes} b$. That is the rows have simply been interchanged.\\
Hence, there exists a $4n_{1}n_{2}\times 4n_{1}n_{2}$ permutation matrix $P$ such that $P \begin{bmatrix}
(I_{n_{1}n_{2}}, 0_{n_{1}n_{2}})^{t} &    0_{2n_{1}n_{2},n_{1}n_{2}}      \\
    0_{2n_{1}n_{2},n_{1}n_{2}}        &  (I_{n_{1}n_{2}}, 0_{n_{1}n_{2}})^{t}
\end{bmatrix} = (I_{2n_{1}n_{2}}, 0_{2n_{1}n_{2}})^{t}$. \\
$P$ being a permutation matrix is invertible and its inverse is $P^{t}$. \\
Now that we know $P$ exists and is invertible, we need to check if the pair of matrices $(P,P)$ is
a map from $(e\widetilde{\otimes} a) \widetilde{\otimes} b$ to $e\widetilde{\otimes} (a \widetilde{\otimes} b)$ in $\mathcal{T}$ and if the pair of matrices $(P^{-1},P^{-1})$ is
a map from $e\widetilde{\otimes} (a \widetilde{\otimes} b)$ to $(e\widetilde{\otimes} a) \widetilde{\otimes} b$. \\
We do it for $(P,P)$, the case of $(P^{-1},P^{-1})$ is completely similar. \\
Once more, $a$ and $b$ are objects of $\mathcal{T}$. So, we can let $a=e^{p}$ and $b=e^{m}$. Observe that $e\widetilde{\otimes} (a \widetilde{\otimes} b)=(e\widetilde{\otimes} a) \widetilde{\otimes} b = (e\widetilde{\otimes} e^{p}) \widetilde{\otimes} e^{m}=e^{1+p+m}=(I_{2^{1+p+m-1}}, I_{2^{1+p+m-1}})= (I_{2^{p+m}}, I_{2^{p+m}})$. Let $r=2^{p+m}.$\\
All we need check now to conclude that $(P,P)$ is a map in $\mathcal{T}$ is that the following diagram commutes:

$$\xymatrix@ R=0.6in @ C=.75in{K[[x]]^{r} \ar[r]^{I_{r}} \ar[d]_{P} &
K[[x]]^{r} \ar[d]^{P} \ar[r]^{I_{r}} & K[[x]]^{r}\ar[d]^{P\,\,\,\,\,\,\,\,\,\,\,\,\,\,\,\,\,\,(\star\star)}\\
K[[x]]^{r} \ar[r]^{I_{r}} & K[[x]]^{r}\ar[r]^{I_{r}} & K[[x]]^{r}}$$

Now, this diagram clearly commutes, so we can take $\zeta:=(P,P)$ and $\zeta^{-1}:=(P^{-1},P^{-1})=(P^{t},P^{t})$.\\
Therefore there exists an isomorphism namely $\zeta$ such that diagram $(2)$ commutes.
\\
\textit{A small remark}: The foregoing proof for the commutativity of diagram $(2)$ helps understand the motivation behind the choice of the objects of $\mathcal{T}$.\\ In fact, if objects were chosen arbitrarily say pairs of matrices $(M,M^{-1})$, as we showed in remark \ref{remark: MF(1) is not semi-unital}, the twin diagram $(\star\star)$ above will commute only if $PM=MP$. But as explained in remark \ref{remark: MF(1) is not semi-unital}, this is not possible as on the left side of the equality, the rows of $M$ are permuted and on the right side the columns are permuted, since $P$ is a permutation matrix. \\
Moreover, though diagram $(\star\star)$ commutes even if $P$ is replaced with any matrix, what we need is a matrix that will make diagram $(2)$ commute and that matrix should also be invertible because we need an isomorphism in diagram $(2)$. \\

The commutativity of diagram $(3)$ is proved in a manner similar to the proof given for the commutativity of
diagram $(2)$.\\
So $e$ is a semi-unit in $(\mathcal{T},\widetilde{\otimes})$.\\
Conclusion: $(\mathcal{T},\widetilde{\otimes})$ is a one-step connected semi-unital semi-monoidal subcategory of $MF(1)$.
\end{proof}
 The above proof works well for $\mathcal{T}$ because the objects of $\mathcal{T}$ are judiciously chosen so that the pair of matrices that make an object in $\mathcal{T}$ is not any kind of matrix and its inverse (in order for the product to yield $1$ times the identity matrix of the right size), but they are identity matrices thanks to which diagrams will be commutative. In fact, diagrams $(2)$ and $(3)$ in definition \ref{defn semi-unital semimonoidal catg}, commute when $a$ and $b$ are objects in $\mathcal{T}$, i.e., of the form $e^{n}$ for some $n\in \mathbb{N}-\{0\}$. But we will see in remark \ref{remark: MF(1) is not semi-unital}, that for arbitrary values of $a$ and $b$, diagrams $(2)$ and $(3)$ are not always commutative. This implies that $(MF(1), \widetilde{\otimes})$ is not a semi-unital semi-monoidal category.
 \begin{remark} \label{remark: MF(1) is not semi-unital}
 We now explain why $(MF(1), \widetilde{\otimes})$ is not a semi-unital semi-monoidal category.
 \begin{enumerate}
   \item We explain that for $a\,\in\, Ob(MF(1))$ of size $n_{1}$, the only reasonable (nonzero) possible choice for $\gamma_{a}: a \rightarrow e \widetilde{\otimes} a$ is what we made for the subcategory $\mathcal{T}$, namely $\gamma_{a}=(\delta',\beta')=((I_{n_{1}},0)^{t}, (I_{n_{1}},0)^{t})$.\\
       First of all, observe that considering the definition of morphisms in $\mathcal{T}$ (i.e., pairs of $(1,0)-$matrices s.t. each column and each row has at most one nonzero entry), the only possible choice for $\gamma_{a}$ in $\mathcal{T}$ is the one we made above (cf. theorem \ref{Thm: T is a semi-unital semi-monoidal catg}), i.e., $\gamma_{a}=((I_{n_{1}},0)^{t}, (I_{n_{1}},0)^{t})$.\\
       It is clear that the only candidate to be a semi-unit in $\mathcal{T}$ was $e=(1,1)$. Hence, it is also the only candidate for $(MF(1), \widetilde{\otimes})$ to be semi-unital. This entails that for $a$ in $Ob(MF(1))$, the only possible way to define $\gamma_{a}$ is $\gamma_{a}=((I_{n_{1}},0)^{t}, (I_{n_{1}},0)^{t})$. Otherwise, $e=(1,1)$ would no more be a semi-unit in $\mathcal{T}$.
\item Next, we prove that with this choice of $\gamma_{a}=((I_{n_{1}},0)^{t}, (I_{n_{1}},0)^{t})$, the diagram $(2)$ above does not commute in general (i.e., for arbitrary values of $a$ and $b$ in $Ob(MF(1))$). That is,
\[
\xymatrix{
a\widetilde{\otimes} b\ar[rr]^{\gamma_{a}\widetilde{\otimes} b} \ar[rd]_{\gamma_{a\widetilde{\otimes} b}} &&
(e\widetilde{\otimes} a) \widetilde{\otimes} b \ar[ld]^{\cong \,\,\,\,\,\,\,\,\,\,\,\,\,\,\,\,(2)}\\
& e\widetilde{\otimes} (a \widetilde{\otimes} b)}
\]
does not commute. In fact, we showed in the proof of theorem \ref{Thm: T is a semi-unital semi-monoidal catg} that with $\gamma=((I_{n_{1}},0)^{t}, (I_{n_{1}},0)^{t})$, the matrix constituting the map $\gamma_{a\widetilde{\otimes} b}$ is permutation equivalent to the matrix constituting the map $\gamma_{a}\widetilde{\otimes} b$. Hence, in order to find the desired isomorphism of diagram $(2)$, all we need do is to find a permutation matrix as explained in the proof of theorem \ref{Thm: T is a semi-unital semi-monoidal catg}. Now, the catch is that we need to verify that this permutation matrix is actually the matrix of a map $ e\widetilde{\otimes} (a \widetilde{\otimes} b) \rightarrow (e\widetilde{\otimes} a) \widetilde{\otimes} b$ in $MF(1)$. It turns out that it is not.\\
Suppose we have already found the permutation matrix that enables us to move from the matrix of $\gamma_{a\widetilde{\otimes} b}$ to the matrix of $\gamma_{a}\widetilde{\otimes} b$, call it $P'$. Now by definition of $\widetilde{\otimes}$, we have $ e\widetilde{\otimes} (a \widetilde{\otimes} b)= (e\widetilde{\otimes} a) \widetilde{\otimes} b$ which is an object of $MF(1)$, so there is a matrix $M$ such that $ e\widetilde{\otimes} (a \widetilde{\otimes} b)= (e\widetilde{\otimes} a) \widetilde{\otimes} b=(M,M^{-1})$. Our aim is to show that $(P',P'): e\widetilde{\otimes} (a \widetilde{\otimes} b) \rightarrow (e\widetilde{\otimes} a) \widetilde{\otimes} b$ is not a map in $MF(1)$ for arbitrary values of $a$ and $b$, because the following diagram cannot commute:

        $$\xymatrix@ R=0.6in @ C=.75in{K[[x]]^{n_{1}} \ar[r]^{M^{-1}} \ar[d]_{P'} &
K[[x]]^{n_{1}} \ar[d]^{P'} \ar[r]^{M} & K[[x]]^{n_{1}}\ar[d]^{P'}\\
K[[x]]^{2n_{1}} \ar[r]^{M^{-1}} & K[[x]]^{2n_{1}}\ar[r]^{M} & K[[x]]^{2n_{1}}}$$
For this diagram to commute, we need to have (from the second square) $P'M=MP'$. Now,
we know that $P'M$ is the matrix obtained from $M$ by permuting the rows according to the permutation $P'$ and $MP'$ is the matrix obtained from $M$ by permuting the columns according to the permutation $P'$. So, $P'M=MP'$ will be true just in case $M$ is the identity matrix. Now, $M$ is not necessarily the identity matrix, for instance if we take $a=(\begin{bmatrix}
4 &   3      \\
1  &  1
\end{bmatrix},\begin{bmatrix}
1 &  -3     \\
-1 &  4
\end{bmatrix})$ and $b=(\begin{bmatrix}
1 &   0      \\
0  &  1
\end{bmatrix},\begin{bmatrix}
1 &  0     \\
0 &  1
\end{bmatrix})$
 \end{enumerate}
 then $M=(1\widetilde{\otimes}\begin{bmatrix}
4 &   3      \\
1  &  1
\end{bmatrix} )\widetilde{\otimes} \begin{bmatrix}
1 &   0      \\
0  &  1
\end{bmatrix} $ is clearly not equal to the identity matrix.
 \end{remark}
\begin{remark} \label{remark: T is not monoidal}
$(\mathcal{T},\widetilde{\otimes})$ is not a monoidal category because it has no unit. In fact, the only candidate to be a unit is $e$. Now, in order to be a unit, $e$ needs to first of all be a pseudo-idempotent (cf. remark \ref{remark on units}). But $e=(1,1)$ is not even a pseudo-idempotent.
We have $e^{2}=(1,1)\widetilde{\otimes} (1,1)=(I_{2}, I_{2})$.\\
Let $\zeta_{1}=(\delta_{1},\beta_{1}): e\rightarrow e^{2}$ and $\zeta_{2}=(\delta_{2},\beta_{2}): e^{2}\rightarrow e$. Consider the following situation:
\newpage

$$\xymatrix@ R=0.6in @ C=.75in{K[[x]] \ar[r]^{1} \ar[d]_{\delta_{1}} &
K[[x]] \ar[d]^{\beta_{1}} \ar[r]^{1} & K[[x]]\ar[d]^{\delta_{1}}\\
K[[x]]^{2} \ar[r]^{I_{2}} & K[[x]]^{2}\ar[r]^{I_{2}} & K[[x]]^{2}\\
}$$
$$
\xymatrix@ R=0.6in @ C=.82in{ \ar[d]_{\delta_{2}} &
 \ar[d]^{\beta_{2}} & \ar[d]^{\delta_{2}}\\
K[[x]] \ar[r]^{1} & K[[x]]\ar[r]^{1} & K[[x]]\\
}$$
From the discussion we had about the choice of matrices constituting $\gamma_{a}$ in remark \ref{remark: MF(1) is not semi-unital}, we have only one (nonzero) choice for $\zeta_{1}$; namely $\zeta_{1}=((1,0)^{t}, (1,0)^{t})$ and similarly we have only one way of defining $\zeta_{2}$; $\zeta_{2}=((1,0), (1,0)): e^{2}\rightarrow e$.
Hence, we clearly obtain $\zeta_{2}\circ\zeta_{1}=((1,0)(1,0)^{t},(1,0)(1,0)^{t})=(1,1)= id_{e}$. Hence,  $\zeta_{2}\circ\zeta_{1} = id_{e}$. \\
But, we do not obtain $\zeta_{1}\circ\zeta_{2}= id_{e^{2}}$. \\
In fact, $\zeta_{1}\circ\zeta_{2}=((1,0)^{t}(1,0),(1,0)^{t}(1,0))=(\begin{bmatrix}
1 &   0      \\
0  &  0
\end{bmatrix} \begin{bmatrix}
1 &   0      \\
0  &  0
\end{bmatrix})\neq (\begin{bmatrix}
1 &   0      \\
0  &  1
\end{bmatrix} \begin{bmatrix}
1 &   0      \\
0  &  1
\end{bmatrix})=id_{e^{2}}$

\end{remark}
Therefore, there is no isomorphism between $e$ and $e^{2}$.
A consequence of remark \ref{remark: T is not monoidal} is that $(MF(1),\widetilde{\otimes} )$ is not a monoidal category since the only candidate to be a unit, namely $e$ is not even a pseudo-idempotent.

\begin{remark} \label{remark: T is not a right-monoidal category}
Moreover, $(\mathcal{T},\widetilde{\otimes})$ is not a right-monoidal category (cf. definition \ref{defn: right monoidal catg}) because when trying to verify if the axioms of definition \ref{defn: right monoidal catg} hold for $\mathcal{T}$, instead of equalities we obtain maps which are not equal but whose representing matrices are row-permutation equivalent. Let us for example illustrate what we mean with the second axiom (cf. definition \ref{defn: right monoidal catg}):
$$\alpha_{e,a,b}\circ \gamma_{a\widetilde{\otimes} b}= \gamma_{a}\widetilde{\otimes} b\,\,\,\,\,\,\,\,\,\,\,\,\,\,\,\,\cdots \,\,\,\,\,\,\,\,\,(Ax.2)$$ where $a,b$ are in $Ob(\mathcal{T})$, $\alpha$ is the associator and $\gamma$ is the natural transformation defined in the proof of theorem \ref{Thm: T is a semi-unital semi-monoidal catg}. If $a$ and $b$ are respectively of sizes $m$ and $n$, then by definition of $\gamma$, $\gamma_{a\widetilde{\otimes} b}=((I_{2mn},0)^{t}, (I_{2mn},0)^{t})$ and since $(e \widetilde{\otimes} a)\widetilde{\otimes} b= e \widetilde{\otimes}(a \widetilde{\otimes} b)$, $\alpha_{e,a,b}= (I_{4mn}, I_{4mn})$ and so the left hand side of $(Ax.2)$ becomes $\alpha_{e,a,b}\circ \gamma_{a\widetilde{\otimes} b}= (I_{4mn}(I_{2mn},0)^{t}, I_{4mn}(I_{2mn},0)^{t})= ((I_{2mn},0)^{t}, (I_{2mn},0)^{t})\,\,\,\,\,\,\,\, \cdots (i)$.\\
Next, by definition \ref{defn zeta tensor tilda X'}, we compute the right hand side of $(Ax.2)$ as follows:
$\gamma_{a}\widetilde{\otimes} b=((I_{m},0)^{t}, (I_{m},0)^{t}) \widetilde{\otimes} b= (\begin{bmatrix}
(I_{m},0)^{t}\otimes I_{n} &   0      \\
0  &  (I_{m},0)^{t}\otimes I_{n}
\end{bmatrix} \begin{bmatrix}
(I_{m},0)^{t}\otimes I_{n} &   0      \\
0  &  (I_{m},0)^{t}\otimes I_{n}
\end{bmatrix})$ \\$= (\begin{bmatrix}
(I_{mn},0)^{t} &   0      \\
0  &  (I_{mn},0)^{t}
\end{bmatrix} \begin{bmatrix}
(I_{mn},0)^{t} &   0      \\
0  &  (I_{mn},0)^{t}
\end{bmatrix})\,\,\,\,\,\,\, \cdots (ii)$\\

The matrices we obtained in $(i)$ and $(ii)$ are row-permutation equivalent but not equal. This proves that $(Ax.2)$ does not hold in $(\mathcal{T}, \widetilde{\otimes})$, so it is not a right-monoidal category. A direct consequence of this result is that $(MF(1), \widetilde{\otimes})$ is not also a right-monoidal category.

\end{remark}

Nevertheless, $(MF(1),\widetilde{\otimes} )$ is still a category which is close to being a monoidal category as we shall see (cf. subsection \ref{subsec: MF(1) is a right pseudo-monoidal catg}).

\subsection{Another application of $\widetilde{\otimes}$: The category $(MF(1), \widetilde{\otimes})$ is a right pseudo-monoidal category} \label{subsec: MF(1) is a right pseudo-monoidal catg}
In this section, we first define what a \textit{right pseudo-monoidal category} is. We observe that this notion is a generalization of the notion of monoidal category. We exploit the results obtained in the previous sections of this paper to show that the category $MF(1)$ is a \textit{right pseudo-monoidal category}. \\

First recall that a semi-monoidal category definition \ref{defn semi-monoidal category} is one endowed with a bifunctor and a natural isomorphism (called the associator \cite{carqueville2016adjunctions}) such that the pentagon diagram (cf. definition \ref{defn semi-monoidal category}) commutes.

\begin{definition} \label{defn right pseudo-monoidal category}
  A \textbf{right pseudo-monoidal category} $\mathcal{C}=<\mathcal{C},\square, e, \alpha, \lambda, \rho>$ is a category $\mathcal{C}$ which possesses a distinguished element $e$, a natural isomorphism $\alpha$ and two natural retractions $\lambda\,and\,\rho$ s.t. the following hold:\\
 $\bullet$ There exists a morphism $\zeta: e^{2}\rightarrow e$ s.t. $\zeta$ has a right inverse.\\
 $\bullet$ $\alpha$ is a natural isomorphism with components $\alpha_{a,b,c}:(a\otimes b)\otimes c \rightarrow a\otimes (b\otimes c)$ such that the following pentagonal diagram

$\xymatrix @ R=0.4in @ C=0.80in{
a\square (b\square (c\square d)) \ar[rr]^{1_{a}\square \alpha}\ar[d]_\alpha &&
a\square ((b\square c)\square d) \ar[d]^\alpha &\\
(a\square b)\square (c\square d) \ar[rd]_\alpha &&
(a\square(b\square c))\square d \ar[ld]^{\alpha \square 1_{d}} &\\
&((a\square b)\square c)\square d} $

commutes for all $a,b,c,d\,\in \,\mathcal{C}$.\\
$\bullet$ $\lambda: e\square (-) \rightarrow (-)$, $\rho: (-)\square e \rightarrow (-)$ are natural\footnote{On p.10 of \cite{carqueville2016adjunctions}, $\lambda$ and $\rho$ are also called left and right unit actions (or unitors). The difference here is that in the definition of a right pseudo-monoidal category, we do not require these unitors to be natural isomorphisms but it is enough for them to have right inverses.} transformations.
\\
$\bullet$ For all objects $a\,\in \mathcal{C}$, $\lambda_{a}: e\square a \rightarrow a$, $\rho_{a}: a\square e \rightarrow a$ have right inverses but do not necessarily have left inverses.
\\
$\bullet$ $\lambda_{e}=\rho_{e}:e\square e \rightarrow e$.\\
$\bullet$ For $a=e$ and for any object $b\,\in \,\mathcal{C}$, the triangular diagram
\[\xymatrix{a\square (e\square b)\ar[rr]^-\alpha \ar[rd]_{1_{a}\square \lambda} &&
(a\square e)\square b \ar[ld]^{\rho\square 1_{b}}\\
&a\square b}\] commutes.

\end{definition}
\begin{remark}
  It is easy to see that every monoidal category (cf. definition \ref{defn monoidal category}) is a right pseudo-monoidal category. In fact, in the foregoing definition, if the triangular diagram commutes for all $a\in \mathcal{C}$; and the maps $\lambda$ and $\rho$ are invertible, then we will recover the definition of a monoidal category. This shows that this notion is a generalization of the classical notion of monoidal category.
\end{remark}
\begin{theorem} \label{thm MF(1) is a pseudo-monoidal catg}
The category $(MF(1), \widetilde{\otimes})$ is a right pseudo-monoidal category.
\end{theorem}
\begin{proof}
Following definition \ref{defn right pseudo-monoidal category}, we need to first of all show that $(MF(1),\widetilde{\otimes})$ is semi-monoidal (cf. definition \ref{defn semi-monoidal category}). Thus we need to show that $\widetilde{\otimes}$ is a bifunctor, and the associator "$\alpha$" in $(MF(1), \widetilde{\otimes})$ is a natural isomorphism such that the pentagon (cf. definition \ref{defn right pseudo-monoidal category}) diagram commutes. \\
Recall that an object of $MF(1)$ is of the form $(M,N)$ where $M=N^{-1}$.\\
In the entire proof; $a,b\,and\,c$ stand for arbitrary objects of $MF(1)$, say $a=e^{p}=(\phi,\psi)$, $b=e^{m}=(\phi',\psi')$ and $c=e^{r}=(\phi'',\psi'')$. \\
$\bullet$ We know that $\widetilde{\otimes}$ is a bifunctor (cf. theorem \ref{thm: mult tens prdt is a bifunctor}).\\
$\bullet$
 It is easy to see that $\alpha=\alpha_{a,b,c}: a\widetilde{\otimes} (b\widetilde{\otimes} c) \xrightarrow[]{=} (a\widetilde{\otimes} b)\widetilde{\otimes} c$ is an identity map and hence it is an isomorphism. It is also easy to see that $\alpha$ is natural for all $a,b,c\,\in \,MF(1)$ and that the above pentagonal diagram commutes, in fact; we actually already proved it above when proving that $\mathcal{T}$ was a semi-monoidal category (cf. theorem \ref{Thm: T is a semi-unital semi-monoidal catg}) . \\

This shows that $(MF(1), \widetilde{\otimes})$ is a semi-monoidal category.\\

Next, we find the distinguished object "$e$" and the morphism $\zeta: e^{2}\rightarrow e$ s.t. $\zeta$ has a right inverse. \\
$\bullet$ Take $e=(1,1)$, the pair of $1\times 1$ matrix factorization. We have $e^{2}=(1,1)\widetilde{\otimes} (1,1)=(I_{2}, I_{2})$.
Consider the following situation:
$$\xymatrix@ R=0.6in @ C=.75in{K[[x]] \ar[r]^{1} \ar[d]_{\delta_{1}} &
K[[x]] \ar[d]^{\beta_{1}} \ar[r]^{1} & K[[x]]\ar[d]^{\delta_{1}}\\
K[[x]]^{2} \ar[r]^{I_{2}} & K[[x]]^{2}\ar[r]^{I_{2}} & K[[x]]^{2}\\
}$$
$$
\xymatrix@ R=0.6in @ C=.82in{ \ar[d]_{\delta_{2}} &
 \ar[d]^{\beta_{2}} & \ar[d]^{\delta_{2}}\\
K[[x]] \ar[r]^{1} & K[[x]]\ar[r]^{1} & K[[x]]\\
}$$

In order to find the map $\zeta$ and its right inverse, it suffices to take: $\zeta=(\delta_{2}=(1,0), \beta_{2}=(1,0)): e^{2}\rightarrow e$, let $\zeta'=(\delta_{1}=(1,0)^{t}, \beta_{1}=(1,0)^{t}): e\rightarrow e^{2}$. Thus, $\zeta\circ\zeta': e \rightarrow e$. Hence, we clearly obtain $\zeta\circ\zeta'=((1,0)(1,0)^{t},(1,0)(1,0)^{t})=(1,1)= id_{e}$ which proves that $\zeta'$ is a right inverse to $\zeta$. \\
$\bullet$ We now show that the maps $\lambda$ and $\rho$ should be natural retractions satisfying $\lambda_{e}=\rho_{e}$. That is, for each $a\,in\, Ob(MF(1))$, $\lambda_{a}$ and $\rho_{a}$ have right inverses and $\lambda_{e}=\rho_{e}$.\\
$\lambda$ is a natural transformation:\\
$\lambda : F=e\widetilde{\otimes} (-) \rightarrow (-)=G$ where G is the identity endofunctor on $MF(1)$ and $F$ is an endofunctor\footnote{It is easy to verify that $F$ is a functor.} on $MF(1)$, such that $F(a)=e\widetilde{\otimes} a$.\\
The family of morphisms $\lambda$ should satisfy the following two requirements:
\begin{enumerate}
  \item For each $a\,in\, Ob(MF(1))$, $\lambda_{a}$ should be a morphism between objects in $MF(1)$.
  Before we proceed, observe that, for any $a=(\phi,\psi)$ of size $n_{1}$ in $MF(1)$,
  $$
(1,1)\widetilde{\otimes}(\phi,\psi)=(
\begin{bmatrix}
   1 \otimes \phi   &          0      \\
    0                  & 1 \otimes \phi
\end{bmatrix},
\begin{bmatrix}
1 \otimes \psi  &    0      \\
    0                 &  1 \otimes \psi
\end{bmatrix})=(\begin{bmatrix}
    \phi  &   0      \\
     0    &  \phi
\end{bmatrix},
\begin{bmatrix}
   \psi  &    0      \\
    0    &  \psi
\end{bmatrix})
$$
  To show that $\lambda_{a}: e\widetilde{\otimes} a\rightarrow a$ is a morphism, we need to find a pair of matrices $(\delta,\beta)$ such that the following diagram commutes:\\
      $$\xymatrix@ R=0.6in @ C=.75in{K[[x]]^{2n_{1}} \ar[r]^{\begin{bmatrix}
\psi  &    0      \\
    0   &  \psi
\end{bmatrix}} \ar[d]_{\delta} &
K[[x]]^{2n_{1}} \ar[d]^{\beta} \ar[r]^{\begin{bmatrix}
\phi  &    0      \\
    0   &  \phi
\end{bmatrix}} & K[[x]]^{2n_{1}}\ar[d]^{\delta\,\,\,\,\,\,\,\,\,\,\,\,\,\,\,\star'}\\
K[[x]]^{n_{1}} \ar[r]^{\psi} & K[[x]]^{n_{1}}\ar[r]^{\phi} & K[[x]]^{n_{1}}}$$
That is,
$$\ast\begin{cases}
 \delta\begin{bmatrix}
\phi  &    0      \\
   0  &  \phi
\end{bmatrix} = \phi\beta  \\
 \psi\delta = \beta\begin{bmatrix}
\psi  &    0      \\
    0   &  \psi
\end{bmatrix}
\end{cases}$$
For $\delta=\beta=(I_{n_{1}},0)$, where $0$ is the zero $n_{1}\times n_{1}$ matrix, the equational system $\ast$ becomes
$$\begin{cases}
 (I_{n_{1}},0)\begin{bmatrix}
\phi  &    0      \\
   0  &  \phi
\end{bmatrix} = \phi (I_{n_{1}},0)  \\
 \psi (I_{n_{1}},0) = (I_{n_{1}},0)\begin{bmatrix}
\psi  &    0      \\
    0   &  \psi
\end{bmatrix}
\end{cases}$$
That is;
$$\begin{cases}
(\phi,0)=(\phi,0) \\
(\psi,0)=(\psi,0)
\end{cases}$$
and this is clearly true. Therefore, we have found a pair of matrices $(\delta,\beta)$ such that diagram $\star'$ commutes, and this means that $\lambda_{a}$  is a map of matrix factorizations.\\
\item Naturality of $\lambda$: \\
Let $b=(\phi',\psi')$ be a matrix factorization of size $n_{2}$ and let $\nu= (\alpha_{\nu}, \beta_{\nu}): a\rightarrow b$ be a map of matrix factorizations. It is easy\footnote{By drawing the twin diagram that has to commute with $(\alpha_{\nu}, \beta_{\nu})$, we see the sizes of $\alpha_{\nu}$ and $\beta_{\nu}$. } to see that $\alpha_{\nu}$ and $\beta_{\nu}$ are each of size $n_{2}\times n_{1}$.
  The following diagram should commute:\\
  \[
\xymatrix{
e\widetilde{\otimes} a\ar[r]^{\lambda_{a}} \ar[d]_{e\widetilde{\otimes} \nu} & a \ar[d]_{\nu}\\
e\widetilde{\otimes} b\ar[r]_{\lambda_{b}} & b}
\]
  \\
  i.e., $\nu \circ \lambda_{a}=\lambda_{b} \circ e\widetilde{\otimes} \nu$ $\cdots (E)$\\
  We know that $e\widetilde{\otimes} a$ is of size $2n_{1}$ since $a$ is of size $n_{1}$. We also know that $\lambda_{a}=[(I_{n_{1}},0),(I_{n_{1}},0)]$. Now by definition of composition of two morphisms in $MF(1)$, the left hand side of equality $(E)$ becomes:\\
  $\nu \circ \lambda_{a}=(\alpha_{\nu}, \beta_{\nu})\circ [(I_{n_{1}},0),(I_{n_{1}},0)]=[\alpha_{\nu}(I_{n_{1}},0),\beta_{\nu}(I_{n_{1}},0)]=[(\alpha_{\nu},0),(\beta_{\nu},0)]$$\cdots \natural$\\
  $0$ in $[(\alpha_{\nu},0),(\beta_{\nu},0)]$ is the $n_{2}\times n_{1}$ zero matrix. \\
  As for the right hand side of $(E)$, first recall that $\lambda_{b}=[(I_{n_{2}},0),(I_{n_{2}},0)]$, (where $0$ is the zero $n_{2}\times n_{2}$ matrix) and by definition \ref{defn zeta is a funtor (right)} of the multiplicative tensor product, we know that $e\widetilde{\otimes} \nu=(1,1)\widetilde{\otimes}(\alpha_{\nu}, \beta_{\nu})=(\begin{bmatrix}
1 \otimes \alpha_{\nu}  &    0      \\
    0                 &  1 \otimes \alpha_{\nu}
\end{bmatrix},
\begin{bmatrix}
1 \otimes \beta_{\nu}  &    0      \\
    0                 &  1 \otimes \beta_{\nu}
\end{bmatrix})=(\begin{bmatrix}
\alpha_{\nu} &    0      \\
    0        &  \alpha_{\nu}
\end{bmatrix},
\begin{bmatrix}
\beta_{\nu}  &    0      \\
    0          &  \beta_{\nu}
\end{bmatrix})$ \\
So, $\lambda_{b} \circ e\widetilde{\otimes} \nu=[(I_{n_{2}},0),(I_{n_{2}},0)]\circ (\begin{bmatrix}
\alpha_{\nu} &    0      \\
    0        &  \alpha_{\nu}
\end{bmatrix},
\begin{bmatrix}
\beta_{\nu}  &    0      \\
    0          &  \beta_{\nu}
\end{bmatrix})\\
=((I_{n_{2}},0)\begin{bmatrix}
\alpha_{\nu} &    0      \\
    0        &  \alpha_{\nu}
\end{bmatrix},(I_{n_{2}},0)\begin{bmatrix}
\beta_{\nu}  &    0      \\
    0          &  \beta_{\nu}
\end{bmatrix})=[(\alpha_{\nu},0),(\beta_{\nu},0)] \cdots \natural\natural$.\\
From $\natural$ and $\natural\natural$, we see that equality $(E)$ holds. That is $\lambda$ is a natural transformation.
\end{enumerate}
$\bullet$ We find the right inverse of $\lambda_{a}$, for any $a=(\phi,\psi)$ of size $n_{1}$ in $MF(1)$.
we denote it $\gamma_{a}: a \rightarrow e\widetilde{\otimes} a$. $\gamma_{a}$ should be a member of the family of morphisms of a natural transformation $\gamma : (-)=G \rightarrow F=e\widetilde{\otimes} (-)$, where G is the identity endofunctor on $MF(1)$ and $F$ is an endofunctor\footnote{this was discussed above when dealing with $\lambda$} on $MF(1)$, such that $F(a)=e\widetilde{\otimes} a$.\\
The family of morphisms $\gamma$ should satisfy the following two requirements:
\begin{enumerate}
  \item For each $a\,in\, Ob(MF(1))$, $\gamma_{a}$ should be a morphism in $MF(1)$.\\
 $\gamma_{a}$ should be a pair of matrices $(\delta',\beta')$ such that the following diagram commutes:\\
      $$\xymatrix@ R=0.6in @ C=.75in{K[[x]]^{n_{1}} \ar[r]^{\psi} \ar[d]_{\delta'} &
K[[x]]^{n_{1}} \ar[d]^{\beta'} \ar[r]^{\phi} & K[[x]]^{n_{1}}\ar[d]^{\delta'\,\,\,\,\,\,\,\,\,\,\,\,\,\,\,\star"}\\
K[[x]]^{2n_{1}} \ar[r]^{\begin{bmatrix}
\psi  &    0      \\
    0   &  \psi
\end{bmatrix}} & K[[x]]^{2n_{1}}\ar[r]^{\begin{bmatrix}
\phi  &    0      \\
    0   &  \phi
\end{bmatrix}} & K[[x]]^{2n_{1}}}$$
That is,
$$\ast'\begin{cases}
 \delta'\phi = \begin{bmatrix}
\phi  &    0      \\
   0  &  \phi
\end{bmatrix}\beta'  \\
 \begin{bmatrix}
\psi  &    0      \\
    0   &  \psi
\end{bmatrix}\delta' = \beta'\psi
\end{cases}$$
For $\delta'=\beta'=(I_{n_{1}},0)^{t}$, where $t$ is the operation of taking the transpose, $0$ is the zero $n_{1}\times n_{1}$ matrix, the equational system $\ast'$ becomes
$$\begin{cases}
 (I_{n_{1}},0)^{t}\phi = \begin{bmatrix}
\phi  &    0      \\
   0  &  \phi
\end{bmatrix} (I_{n_{1}},0)^{t}  \\
 \begin{bmatrix}
\psi  &    0      \\
    0   &  \psi
\end{bmatrix} (I_{n_{1}},0)^{t} = (I_{n_{1}},0)^{t}\psi
\end{cases}$$
That is;
$$\begin{cases}
(\phi,0)^{t}=(\phi,0)^{t} \\
(\psi,0)^{t}=(\psi,0)^{t}
\end{cases}$$
and this is clearly true. Therefore, we have found a pair of matrices $(\delta',\beta')$ such that diagram $\star"$ commutes, and this means that $\gamma_{a}$  is a map of matrix factorizations.
 \item Naturality of $\gamma$:\\
 Let $b=(\phi',\psi')$ be a matrix factorization of size $n_{2}$ and let $\mu= (\alpha_{\mu}, \beta_{\mu}): a\rightarrow b$ be a map of matrix factorizations. It is easy\footnote{By drawing the twin diagram that has to commute with $(\alpha_{\mu}, \beta_{\mu})$, we see the sizes of $\alpha_{\mu}$ and $\beta_{\mu}$. } to see that $\alpha_{\mu}$ and $\beta_{\mu}$ are each of size $n_{2}\times n_{1}$.
  The following diagram should commute:\\
  \[
\xymatrix{
 a\ar[r]^{\gamma_{a}} \ar[d]_{\mu} & e\widetilde{\otimes} a \ar[d]_{e\widetilde{\otimes}\mu}\\
 b\ar[r]_{\gamma_{b}} & e\widetilde{\otimes} b}
\]
  \\
  i.e., $e\widetilde{\otimes}\mu \circ \gamma_{a}=\gamma_{b} \circ \mu$ $\cdots (E')$\\
  We know that $e\widetilde{\otimes} a$ is of size $2n_{1}$ since $a$ is of size $n_{1}$. We also know that $\gamma_{b}=[(I_{n_{2}},0)^{t},(I_{n_{2}},0)^{t}]$. Now by definition of composition of two morphisms in $MF(1)$, the right hand side of equality $(E')$ becomes:\\
  $\gamma_{b}\circ \mu= [(I_{n_{2}},0)^{t},(I_{n_{2}},0)^{t^{}}]\circ (\alpha_{\mu}, \beta_{\mu})=[(I_{n_{2}},0)^{t}\alpha_{\mu}, (I_{n_{2}},0)^{t}\beta_{\mu}]=[(\alpha_{\mu},0)^{t},(\beta_{\mu},0)^{t}]$$\cdots \natural'$\\
  $0$ in $[(\alpha_{\mu},0)^{t},(\beta_{\mu},0)^{t}]$ is the $n_{2}\times n_{1}$ zero matrix. \\
  As for the left hand side of $(E')$, first recall that $\gamma_{a}=[(I_{n_{1}},0)^{t},(I_{n_{1}},0)^{t}]$, (where $0$ is the zero $n_{1}\times n_{1}$ matrix) and by definition \ref{defn zeta is a funtor (right)} of the multiplicative tensor product, we know that $e\widetilde{\otimes} \mu=(1,1)\widetilde{\otimes}(\alpha_{\mu}, \beta_{\mu})=(\begin{bmatrix}
1 \otimes \alpha_{\mu}  &    0      \\
    0                 &  1 \otimes \alpha_{\mu}
\end{bmatrix},
\begin{bmatrix}
1 \otimes \beta_{\mu}  &    0      \\
    0                 &  1 \otimes \beta_{\mu}
\end{bmatrix})=(\begin{bmatrix}
\alpha_{\mu} &    0      \\
    0        &  \alpha_{\mu}
\end{bmatrix},
\begin{bmatrix}
\beta_{\mu}  &    0      \\
    0          &  \beta_{\mu}
\end{bmatrix})$ \\
So, $e\widetilde{\otimes} \mu \circ \gamma_{a} = (\begin{bmatrix}
\alpha_{\mu} &    0      \\
    0        &  \alpha_{\mu}
\end{bmatrix},
\begin{bmatrix}
\beta_{\mu}  &    0      \\
    0          &  \beta_{\mu}
\end{bmatrix})\circ [(I_{n_{1}},0)^{t},(I_{n_{1}},0)^{t}]\\
=(\begin{bmatrix}
\alpha_{\mu} &    0      \\
    0        &  \alpha_{\mu}
\end{bmatrix}(I_{n_{1}},0)^{t},\begin{bmatrix}
\beta_{\mu}  &    0      \\
    0          &  \beta_{\mu}
\end{bmatrix}(I_{n_{1}},0)^{t})=[(\alpha_{\mu},0)^{t},(\beta_{\mu},0)^{t}] \cdots \natural\natural'$.\\
From $\natural'$ and $\natural\natural'$, we see that equality $(E')$ holds. That is $\gamma$ is a natural transformation.
\end{enumerate}
Next, we show that $\gamma_{a}$ is the right inverse of $\lambda_{a}$ by computing the following: $\lambda_{a}\circ \gamma_{a}=[(I_{n_{1}},0),(I_{n_{1}},0)]\circ [(I_{n_{1}},0)^{t},(I_{n_{1}},0)^{t}]= (I_{n_{1}},I_{n_{1}})=id_{a}$. So $\gamma_{a}$ is the right inverse of $\lambda_{a}$.\\\\

$\bullet$ To see that $\rho$ is a natural transformation and that for any objet $a$ in $MF(1)$, $\rho_{a}$ has a right inverse, it suffices to observe that both $\lambda_{a}$ and $\rho_{a}$ have the same domain and codomain since for any $a=(\phi,\psi)$ in $MF(1)$, we have: \\
$$
(\phi,\psi)\widetilde{\otimes} (1,1)=(
\begin{bmatrix}
    \phi\otimes 1  &          0      \\
    0                  &\phi\otimes 1
\end{bmatrix},
\begin{bmatrix}
   \psi\otimes 1  &    0      \\
    0                 &  \psi\otimes 1
\end{bmatrix})=(\begin{bmatrix}
    \phi  &   0      \\
     0    &  \phi
\end{bmatrix},
\begin{bmatrix}
   \psi  &    0      \\
    0    &  \psi
\end{bmatrix})
$$
Similarly,
$$
(1,1)\widetilde{\otimes}(\phi,\psi)=(
\begin{bmatrix}
   1 \otimes \phi   &          0      \\
    0                  & 1 \otimes \phi
\end{bmatrix},
\begin{bmatrix}
1 \otimes \psi  &    0      \\
    0                 &  1 \otimes \psi
\end{bmatrix})=(\begin{bmatrix}
    \phi  &   0      \\
     0    &  \phi
\end{bmatrix},
\begin{bmatrix}
   \psi  &    0      \\
    0    &  \psi
\end{bmatrix})
$$

So, we define $\rho_{a}=\lambda_{a}$ for any $a$ in $MF(1)$.\\

We also clearly have $\rho_{e}=\lambda_{e}$.\\\\
$\bullet$ Finally, for any object $b\,\in \, MF(1)$ and for $a=e$, we prove that the following triangular diagram commutes:
\[
\xymatrix{
a\widetilde{\otimes} (e\widetilde{\otimes} b)\ar[rr]^-\alpha_{a,e,b} \ar[rd]_{1_{a}\widetilde{\otimes} \lambda_{b}} &&
(a\widetilde{\otimes} e)\widetilde{\otimes} b \ar[ld]^{\rho_{a}\widetilde{\otimes} 1_{b}\,\,\,\,\,\,\,\,\,\,\,\,\,(\star''')}\\
&a\widetilde{\otimes} b}
\]
Our goal here is to show that the diagram $(\star''')$ commutes i.e., $\rho_{a}\widetilde{\otimes} 1_{b}\circ \alpha_{a,e,b} = 1_{a}\widetilde{\otimes}\lambda_{b}$ i.e., $\rho_{a}\widetilde{\otimes} 1_{b}= 1_{a}\widetilde{\otimes}\lambda_{b}$ since the associator $\alpha$ is the identity. \\
We use definition \ref{defn zeta is a bifuntor} to verify that this equality holds.\\

$\rho_{a}\widetilde{\otimes} 1_{b}= (\begin{bmatrix}
(I_{n_{1}},0)\otimes I_{n_{2}}  &    0      \\
    0                 &  (I_{n_{1}},0)\otimes I_{n_{2}}
\end{bmatrix},\begin{bmatrix}
(I_{n_{1}},0)\otimes I_{n_{2}}  &    0      \\
    0                 &  (I_{n_{1}},0)\otimes I_{n_{2}}
\end{bmatrix})\\=(\begin{bmatrix}
(1,0)\otimes I_{n_{2}}  &    0      \\
    0                 &  (1,0)\otimes I_{n_{2}}
\end{bmatrix},\begin{bmatrix}
(1,0)\otimes I_{n_{2}}  &    0      \\
    0                 &  (1,0)\otimes I_{n_{2}}
\end{bmatrix})\, since\,n_{1}=1\,as\,a=e\\= (\begin{bmatrix}
(I_{n_{2}},0)  &    0      \\
    0                 &  (I_{n_{2}},0)
\end{bmatrix},\begin{bmatrix}
(I_{n_{2}},0)  &    0      \\
    0                 &  (I_{n_{2}},0)
\end{bmatrix})$ \,\,\,\,\,\,\,\,\,\,\,\,\,$\cdots \flat$ \\

$1_{a}\widetilde{\otimes}\lambda_{b}=(\begin{bmatrix}
I_{n_{1}}\otimes (I_{n_{2}},0)  &    0      \\
    0                 &  I_{n_{1}}\otimes (I_{n_{2}},0)
\end{bmatrix},\begin{bmatrix}
I_{n_{1}}\otimes (I_{n_{2}},0)  &    0      \\
    0                 &  I_{n_{1}}\otimes (I_{n_{2}},0)
\end{bmatrix})\\= (\begin{bmatrix}
1\otimes (I_{n_{2}},0)  &    0      \\
    0                 &  1\otimes (I_{n_{2}},0)
\end{bmatrix},\begin{bmatrix}
1\otimes (I_{n_{2}},0)  &    0      \\
    0                 &  1\otimes (I_{n_{2}},0)
\end{bmatrix})\, since\,n_{1}=1\,as\,a=e\\= (\begin{bmatrix}
(I_{n_{2}},0)  &    0      \\
    0               &  (I_{n_{2}},0)
\end{bmatrix},\begin{bmatrix}
 (I_{n_{2}},0)  &    0      \\
    0                 &  (I_{n_{2}},0)
\end{bmatrix})$ \,\,\,\,\,\,\,\,\,\,\,\,\,$\cdots \flat'$\\

From $\flat$ and $\flat'$, it is clear that $\rho_{a}\widetilde{\otimes} 1_{b}= 1_{a}\widetilde{\otimes}\lambda_{b}$.\\
Therefore $(MF(1), \widetilde{\otimes})$ is a right pseudo-monoidal category. QED

\end{proof}
\begin{remark}
  When proving the commutativity of the triangular diagram in the foregoing proof, we kept writing $a$ instead of directly writing $e$ because we wanted to point out the fact that this diagram is simply the triangular diagram one has in the definition of a monoidal category, except that here, the diagram commutes only for $a=e$. It is easy to see that if $a\neq e$ (meaning $n_{1}\neq 1$), then $\rho_{a}\widetilde{\otimes} 1_{b}\neq 1_{a}\widetilde{\otimes}\lambda_{b}$. In fact, the pair of matrices representing these two maps $\rho_{a}\widetilde{\otimes} 1_{b}\,and\, 1_{a}\widetilde{\otimes}\lambda_{b}$ will be permutation similar but not equal. So, $(MF(1), \widetilde{\otimes})$ resembles a monoidal category in many respects without being one. That is one of the motivations behind the appellation right pseudo-monoidal category.
\end{remark}


\begin{quote}
  \textbf{Acknowledgments}
\end{quote}
This work was carried out while doing my Ph.D. at the University of Ottawa in Canada. I would like to thank Prof. Dr. Richard Blute who was my Ph.D. supervisor for the fruitful interactions.
This research was supported in part by the Bank of Montreal financial group award I.

\bibliography{fomatati_ref}
\addcontentsline{toc}{section}
{References}

\end{document}